# RANDOM WALK IN MARKOVIAN ENVIRONMENT

By Dmitry Dolgopyat,[1] Gerhard Keller[2] and Carlangelo Liverani[3]

*University of Maryland, University of Erlangen-Nuremberg and University of Rome Tor Vergata*

We prove a quenched central limit theorem for random walks with bounded increments in a randomly evolving environment on $\mathbb{Z}^d$. We assume that the transition probabilities of the walk depend not too strongly on the environment and that the evolution of the environment is Markovian with strong spatial and temporal mixing properties.

**1. Introduction.** The study of random walks in random environment encompasses a considerable range of possibilities that have been addressed in an extensive body of literature. We refer to [25, 26] for recent reviews of the field. Here, we consider a situation in which the environment is not static, but has an evolution with strong mixing properties and the transition probabilities of the random walk have a weak dependence on the environment. Note, however, that we have an explicit bound on how strong the dependence on the environment may be. In this case the situation is simpler than in the case of static environment; indeed, we will see that the phenomena known as Sinai traps [23] cannot take place.

Random walks in dynamical environment have been intensively studied under various assumptions (see, e.g., [1, 2, 3, 4, 5, 6, 7, 10, 14, 17, 20, 22, 23, 24]). In fact, [7] considers quite general statical environments and even though it does not formally cover dynamical environments, there seems no conceptual difficulty in doing so. Here, we will consider a finite-range walk in

Received February 2007; revised October 2007.

[1]Supported by the NSF and the IPST.

[2]Supported by the DFG.

[3]Supported in part by M.I.U.R. (Prin 2004028108).

$\mathbb{Z}^d$ with the environment being a (rather general) space–time mixing Markov chain. This generalizes the case, well studied in the literature [1, 2, 3, 4], in which the Markov chain has a product structure; that is, at each site of the lattice a time-mixing Markov chain acts independently on the other sites. In the latter situation, it has been proven that the random walk satisfies an almost sure quenched (i.e., where the histories of the environment are held fixed) CLT for each $d \geq 3$; see [1, 4]. Here, we prove the same result for each $d \geq 1$ and for more general classes of environments.

The general strategy of the proof follows the well-established path of considering the process as seen from the particle and studying such a process via a martingale approximation [11, 16, 18, 19]. In particular, we use the work of [19] in the spirit of [8, 9, 20]. Yet, as we do not discuss the invariance principle and only consider finite-range walks, our arguments are a bit simpler and more direct than those in [20].

A pleasant feature of our approach is that by making heavier use of dynamical arguments, we are able to employ the same methods to establish the mixing properties of the environment and to prove the quenched CLT. In fact, we first prove a CLT under some abstract conditions, then introduce the class of environments and proceed to prove that the above conditions are satisfied.

The paper is organized as follows. In Section 2, we first describe a model in which the process of the environment satisfies a certain number of abstract conditions, and we prove the quenched CLT (Theorem 1), provided a certain correlation decay estimate, (2.21), can be verified. Then, in Section 2.4, we describe a class of models that are claimed (Theorem 2) to satisfy the above abstract conditions. In Section 3, we first show that the inequality (2.21) is equivalent to an estimate for two independent random walks evolving in the same environment (see [5]). The rest of the section is devoted to proving such an estimate. In Section 4, we show that the abstract condition under which the almost sure quenched CLT has been proven before are in fact satisfied by the aforementioned large class of Markov environments, provided the dependence on the environment is sufficiently weak. This proves Theorem 2. Finally, in Appendix, we recall some facts from [19] and slightly generalize some estimates from that paper that we need for our proofs.

Convention. In this paper, we will use $C$ to designate a generic constant depending only on the quantities appearing in the Assumptions (A0)–(A8) below. We will use $C_{a,b,c,\ldots}$ for constants also depending on parameters $a, b, c, \ldots$. Consequently, the actual numerical value of such constants may vary from one occurrence to the next. On the contrary, we will use $C_0, C_1, \ldots,$ to designate constants whose value is held fixed through the paper.

## 2. Model and results.

2.1. *The random walk in random environment.* Let $I$ be a compact Polish space (including the possibility of being finite or countable) and $\theta = (\theta^q)_{q\in\mathbb{Z}^d} \in \Theta := I^{\mathbb{Z}^d}$ be an environment on $\mathbb{Z}^d$. We equip $\Theta$ with the product topology and the Borel $\sigma$-algebra. We assume that the environment has a Markovian time evolution. Let $(\theta_t)_{t\in\mathbb{N}}$ be such a Markov process so that $\theta_t$, with values in $\Theta$, is the environment at time $t \in \mathbb{N}$.

We will use the notation $(\theta_t)^q$, $q \in \mathbb{Z}^d$, with values in $I$, to designate the space components of $\theta_t$ at position $q$. As usual, we will often use the same notation for the random variable and its values since it creates no confusion.

In other words, we have a Markov process with transition probabilities, for each measurable set $A \subset \Theta$,

$$P(\{\theta_{t+1} \in A\} \mid \theta_t) = p(\theta_t, A). \tag{2.1}$$

We require the process to be Feller and translation invariant, that is, $p(\theta, A) = p(\tau^z\theta, \tau^z A)$ for each $\theta, z, A$, where $(\tau^z\theta)^q = \theta^{q+z} \in I$. We will call $\mathbb{P}^e_\nu$ the *measure* on the set $\Omega := \Theta^{\mathbb{N}}$ of environment histories generated by the process (2.1) started with the initial measure $\nu$ on $\Theta$, while we use $\mathbb{P}^e_\theta$ if the process is started in the configuration $\theta \in \Theta$. We will use $\mathbb{E}_{\mathbb{P}^e_\nu}$ for the expectation with respect to $\mathbb{P}^e_\nu$. Note that the translation invariance of the kernel $p$ implies translation equivariance of the measures $\mathbb{P}^e_\theta$, namely $\mathbb{P}^e_{\tau^z\theta}(\tau^z A) = \mathbb{P}^e_\theta(A)$ for $A \subseteq \Omega$ and where $\tau^z$ acts on $\Omega$ by pure space translation.

We then consider a random walk $X_t$ started at $X_0 = 0$ in such an environment. More precisely, let $\Lambda := \{z \in \mathbb{Z}^d : \|z\| \leq C_1\}$ and $\Delta_{t+1} := X_{t+1} - X_t$ (here, and in the following, $\|v\|$ means $\sup_i |v_i|$). The process $(X_t, \theta_t)_{t\in\mathbb{N}}$ is then defined by the transition probabilities

$$P(\{\Delta_{t+1} = z, \theta_{t+1} \in A\} \mid X_t, \theta_t) = \pi_z(\tau^{X_t}\theta_t)p(\theta_t, A), \tag{2.2}$$

where $\pi_z \equiv 0$ for $z \notin \Lambda$, and $\pi_z(\theta)$ depends on $\theta$ only through $(\theta^q)_{q\in\Lambda}$ and is continuous as a function of these variables.

The basic space on which all processes studied in this paper can be defined is $\Omega \times \Lambda^{\mathbb{N}}$, with elements $((\theta_t)_{t\in\mathbb{N}}, (\Delta_t)_{t\in\mathbb{N}})$. The probability measures $\mathbf{P}^e_\nu$ on this space we are interested in are skew product measures with 'base' $\mathbb{P}^e_\nu$ on $\Omega$ and the transition kernel $\mathbb{P}_{(\theta_t)}$, which is the distribution of the increments of the walk on a given space–time environment $(\theta_t)_{t\in\mathbb{N}}$.

It is well known that to study the properties of $X_t$, it is convenient to study the process of the environment as seen from the particle. In fact, such a process can be considered in several ways, two of which will be relevant in the sequel.

2.2. *The process of the environment as seen from the particle.* We look at the environment history $(\theta_s)_{s\in\mathbb{N}}$, not from the origin of the lattice, but from the random position of the particle, and use the letter $\boldsymbol{\omega}$ to denote it. Formally, $\boldsymbol{\omega}$ is also an element of $\Omega$, but the interpretation is different. On $\Omega$, we define the space–time translations $\boldsymbol{\tau}^{a,b}:\Omega\to\Omega$, namely if $\boldsymbol{\omega}=(\theta_t^q)_{t\in\mathbb{N},q\in\mathbb{Z}^d}$ and $(\tilde\theta_t^q)_{t\in\mathbb{N},q\in\mathbb{Z}^d}:=\tau^{a,b}\boldsymbol{\omega}$, then $\tilde\theta_t^q=\theta_{t+a}^{q+b}$.

Let us call $\boldsymbol{\Omega}:=\Omega^{\mathbb{N}}$ the set of all possible paths of space–time histories. $\Omega$ and $\boldsymbol{\Omega}$ are equipped with the obvious product topologies and the corresponding Borel $\sigma$-algebras. As $I$ is separable, these Borel $\sigma$-algebras are at the same time product $\sigma$-algebras, so, for example, the Borel $\sigma$-algebra on $\boldsymbol{\Omega}$ is the product of the one on $\Omega$. In order to describe the process of the environment as seen from the particle, we define the measurable map $\Phi:\Omega\times\Lambda^{\mathbb{N}}\to\boldsymbol{\Omega}$,

$$\Phi((\theta_t)_{t\in\mathbb{N}},(\Delta_t)_{t\in\mathbb{N}})=(\boldsymbol{\omega}_n)_{n\in\mathbb{N}}\qquad\text{with } \boldsymbol{\omega}_0=(\theta_t)_{t\in\mathbb{N}},\boldsymbol{\omega}_n=\boldsymbol{\tau}^{n,X_n}\boldsymbol{\omega}_0. \tag{2.3}$$

It transforms a measure $\mathbf{P}_\nu^e$ into the measure $\mathbf{P}_\nu:=\mathbf{P}_\nu^e\circ\Phi^{-1}$ on $\boldsymbol{\Omega}$. This is the distribution of the process $(\boldsymbol{\omega}_n)_{n\in\mathbb{N}}$ of space–time histories, as seen from the particle under the basic probability measure $\mathbf{P}_\nu^e$. The map $\Phi$ is an almost sure bijection between the probability spaces $(\Omega\times\Lambda^{\mathbb{N}},\mathbf{P}_\nu^e)$ and $(\boldsymbol{\Omega},\mathbf{P}_\nu)$, provided the set of $\boldsymbol{\tau}^{0,b}$-invariant space–time histories $(\theta_t)_{t\in\mathbb{N}}$ has $\mathbf{P}_\nu^e$-measure zero for all $b$. Hence, it is simply a matter of convenience on which basic space we interpret our random variables. For convenience, we also introduce the random variables $\omega_t=(\boldsymbol{\omega}_t)_0$ for each $t\in\mathbb{N}$. Observe that $\omega_t=\tau^{X_t}\theta_t$ are elements of $\Theta$.

In the following lemma, we collect some properties of the processes $(\omega_t)_{t\in\mathbb{N}}$ and $(\boldsymbol{\omega}_t)_{t\in\mathbb{N}}$. The proof is by simple direct computation.

Here, and in the following, we will use $\mathcal{C}^0$ to denote the space of continuous functions and $\mathcal{C}^0_{\text{loc}}$ for the continuous functions depending only on finitely many variables.

LEMMA 2.1. *Let $\nu$ be any initial measure on $\Theta$ and let $\mathbf{P}_\nu$ be the measure on $\boldsymbol{\Omega}$ constructed from it as described above (i.e., via the intermediate steps $\mathbb{P}_\nu^e$ and $\mathbf{P}_\nu^e$).*

(1) $(\omega_t)_{t\in\mathbb{N}}$ *is a Markov process with transition probabilities*

$$\mathbf{P}_\nu(\{\omega_{t+1}\in A\}|\omega_t)=\sum_{z\in\Lambda}\pi_z(\omega_t)p(\omega_t,\tau^{-z}A) \tag{2.4}$$

*and Feller Markov operator $S:\mathcal{C}^0(\Theta)\to\mathcal{C}^0(\Theta)$ defined by*

$$Sf(\omega):=\sum_{z\in\Lambda}\int_\Theta f(\tau^z\omega')\pi_z(\omega)p(\omega,d\omega')=\mathbb{E}_{\mathbf{P}_\nu}(f(\omega_{t+1})|\omega_t=\omega). \tag{2.5}$$

(2) *$(\boldsymbol{\omega}_t)_{t\in\mathbb{N}}$ is a Markov process with transition probabilities*

$$\mathbf{P}_\nu(\{\boldsymbol{\omega}_{t+1}\in A\}|\boldsymbol{\omega}_t)=\sum_{z\in\Lambda}\pi_z((\boldsymbol{\omega}_t)_0)\mathbb{1}_A(\boldsymbol{\tau}^{1,z}\boldsymbol{\omega}_t) \tag{2.6}$$

*and Feller Markov operator $\Pi:\mathcal{C}^0(\Omega)\to\mathcal{C}^0(\Omega)$ defined by*

$$\Pi f(\boldsymbol{\omega}):=\sum_{z\in\Lambda}\pi_z(\theta_0)f(\boldsymbol{\tau}^{1,z}\boldsymbol{\omega}), \tag{2.7}$$

*with the notation $\boldsymbol{\omega}=(\theta_t)_{t\in\mathbb{N}}$ explained above.*

To successfully use both types of processes, the original one and the one seen from the particle, it will be necessary to have initial measures which result in ergodic stationary processes. More precisely, we assume the following.

Assumption (A0) (*Mixing*). There exist unique measures $\mu_e$ and $\mu$ on $\Theta$ such that the processes (2.1) and (2.4), started with the initial distribution $\mu_e$ and $\mu$, respectively, are stationary, ergodic and mixing. In addition, $\mu_e$ is not supported on the translation invariant configurations.

Assumption (A1) (*Absolute continuity*). The measures $\mu$ and $\mu_e$ are equivalent.

In particular, the measure $\mu$ is uniquely characterized by the stationarity condition $\mathbb{E}_\mu(Sf)=\mathbb{E}_\mu(f)$ for all $f\in\mathcal{C}^0(\Theta)$.

Both measures $\mu_e$ and $\mu$ from Assumption (A0) can be used as starting measures for the process $(\theta_t)_{t\in\mathbb{N}}$ of the environment, thus giving rise to measures $\mathbb{P}^e_{\mu_e}$ and $\mathbb{P}^e_\mu$ on $\Omega$.

Clearly, $\mathbb{P}^e_{\mu_e}$ is stationary, and the corresponding measure $\mathbf{P}^e_{\mu_e}$ on $\Omega\times\Lambda^{\mathbb{N}}$, defined in Section 2.1, is our basic reference probability that we will denote simply by $\mathbf{P}^e$. In contrast, the measure $\mathbb{P}^e_\mu$ is not stationary in general, but if we use it to define the measure $\mathbf{P}^e_\mu$ on $\Omega\times\Lambda^{\mathbb{N}}$, then the corresponding measure $\mathbf{P}:=\mathbf{P}^e_\mu\circ\Phi^{-1}$ on $\boldsymbol{\Omega}$ is stationary; in other words, the process $(\boldsymbol{\omega}_t)_{t\in\mathbb{N}}$ has the stationary distribution $\mathbf{P}$ under the probability $\mathbf{P}^e_\mu$. Indeed, a direct computation which uses the translation equivariance of the probability kernel $p$ shows the following.

Lemma 2.2. *Under Assumption (A0), $\mathbb{E}_{\mathbb{P}^e_\mu}(\Pi h)=\mathbb{E}_{\mathbb{P}^e_\mu}(h)$ for each $h\in\mathcal{C}^0(\Omega)$. As $\Pi$ does not increase the supremum-norm, it follows, in particular, that $\Pi$ is an $L^2(\Omega,\mathbb{P}^e_\mu)$-contraction.*

REMARK 2.3. As the measures $\mu_e$ and $\mu$ on $\Theta$ are equivalent, the Markov measures $\mathbb{P}^e_{\mu_e}$ and $\mathbb{P}^e_{\mu}$ on $\Omega$ and the measures $\mathbf{P}^e = \mathbf{P}^e_{\mu_e}$ and $\mathbf{P}^e_{\mu}$ on $\Omega \times \Lambda^{\mathbb{N}}$ are also equivalent. It follows that the same is true for the corresponding measures $\mathbf{P}^e \circ \Phi^{-1}$ and $\mathbf{P} = \mathbf{P}^e_{\mu} \circ \Phi^{-1}$ on $\mathbf{\Omega}$. Therefore, all statements concerning almost sure behavior of our processes have the same meaning, regardless of the measure we are referring to. Observe, however, that this does not mean that the obvious projections of $\mathbf{P}$ and $\mathbf{P}^e$ to $\Omega$ are equivalent.

2.3. *A general quenched CLT.* We assume the following.

ASSUMPTION (A2) (*Time-mixing of the environment as seen from the particle*). There exists $\eta < 1$ such that for each $\varphi \in \mathcal{C}^0_{\mathrm{loc}}(\Theta)$ depending on $M$ variables and for each $n \in \mathbb{N}$,

$$\|S^n \varphi - \mathbb{E}_{\mu}(\varphi)\|_{\infty} \leq CM\eta^n \|\varphi\|_{\infty}.$$

ASSUMPTION (A3) (*Space-mixing of the environment*). There exists $\xi > 4$ such that if $\psi \in \mathcal{C}^0(\Theta)$, $\varphi \in \mathcal{C}^0_{\mathrm{loc}}(\Theta)$ and the supports of $\varphi$ and $\psi$ are at a distance $L$, then

$$|\mu_e(\varphi\psi) - \mu_e(\varphi)\mu_e(\psi)| \leq C_{\varphi} L^{-\xi} \|\varphi\|_{\infty} \|\psi\|_{\infty}.$$

ASSUMPTION (A4) (*Locality of environment dynamics*). There exist $\xi > 4$ and $\tilde{\xi} > 0$ such that for all $M, L, s \in \mathbb{N}$ and $A, B \subset \mathbb{Z}^d$ with diameter at most $M$ and distance $d(A, B) > L$, for all $f, g : \Omega \to \mathbb{R}$ such that $f$ depends only on variables in $A^{\{0,\dots,s\}}$ and $g$ only on variables in $B^{\{0,\dots,s\}}$, and for each $\theta \in \Theta$, the following inequality holds:

$$|\mathbb{E}_{\mathbb{P}^e_\theta}(f(\theta_1, \dots, \theta_s) g(\theta_1, \dots, \theta_s)) - \mathbb{E}_{\mathbb{P}^e_\theta}(f)\mathbb{E}_{\mathbb{P}^e_\theta}(g)| \leq C_M s^{\tilde{\xi}} L^{-\xi} \|f\|_{\infty} \|g\|_{\infty}.$$

ASSUMPTION (A5) (*Ellipticity*). There exist $\gamma_z \geq 0$, $c > 0$ with $\sum_{z\in\Lambda} \gamma_z = 1$ and $|\sum_{z\in\Lambda} \gamma_z e^{i\langle l,z\rangle}| < 1$ for any $l \in \mathbb{Z}^d \setminus \{0\}$, such that $\pi_z(\theta) \geq c\gamma_z$ for $\mathbb{P}^e_{\mu_e}$-almost every $\theta$. In the following, we will set $\gamma = c \inf\{\gamma_z \neq 0\} > 0$.

LEMMA 2.4. *Assumptions (A0), (A1) and (A3) imply that the translation invariant environment configurations have zero $\mathbb{P}^e_{\mu}$-measure.*

PROOF. Let $d$ be a metric on $I$. Next, given $b \in \mathbb{Z}^d$, let $A_{M,\delta} := \{\theta \in \Theta : d((\tau^b\theta)_q, (\theta)_q) \leq \delta \ \forall \|q\| \leq M\}$, $A := \{\theta \in \Theta : \tau^b\theta = \theta\}$. Then, for each $n \in \mathbb{N}$ and $\varepsilon > 0$,

$$\mu_e(A_{M,\delta}) - \varepsilon \leq \mu_e(A) = \mu_e(\mathbb{1}_A \tau^{nb} \mathbb{1}_A) \leq \mu_e(\mathbb{1}_{A_{M,\delta}} \tau^{nb} \mathbb{1}_{A_{M,\delta}}),$$

provided $\delta$ is small enough and $M$ large enough. Choosing $n$ sufficiently large Assumption (A3) implies

$$\mu_e(A) \leq \mu_e(A)^2 + C_{M,\delta} \|bn\|^{-\xi} + C\varepsilon.$$

By the arbitrary nature of $n$ and $\varepsilon$, it follows that $\mu_e(A) \in \{0,1\}$, but $\mu_e(A) = 1$ is ruled out by Assumption (A0), hence, we must have $\mu_e(A) = 0$ and $\mathbb{P}^e_{\mu_e}(\{\theta_t \in A\}) = \mu_e(A) = 0$ for each $t \in \mathbb{N}$. The claim then follows by Assumption (A1). □

REMARK 2.5. Lemma 2.4 implies that the map $\Phi : \Omega \times \Lambda^{\mathbb{N}} \to \mathbf{\Omega}$ is indeed an almost sure bijection. Therefore, the $\sigma$-algebra $\mathcal{F}_t := \sigma\{\boldsymbol{\omega}_0, \dots, \boldsymbol{\omega}_t\}$ and the $\sigma$-algebra $\sigma\{\Delta_1, \dots, \Delta_t, (\theta_s)_{s\in\mathbb{N}}\}$ coincide $\mathbf{P}$-almost surely. Similarly, $\Delta_{t+1}$ is $\sigma(\boldsymbol{\omega}_t, \boldsymbol{\omega}_{t+1})$-measurable.

As $(\boldsymbol{\omega}_t)$ is a Markov process under $\mathbf{P}$, conditional expectations of the form $\mathbb{E}(G|\mathcal{F}_t)$ can be written as functions of $\boldsymbol{\omega}_t$ alone if $G$ is $\sigma(\boldsymbol{\omega}_t, \boldsymbol{\omega}_{t+1}, \dots)$-measurable. This applies, in particular, to $\Delta_{t+1}$. Hence,

$$\mathbb{E}(\Delta_{t+1}|\mathcal{F}_t) = \sum_{z\in\Lambda} z\pi_z(\tau^{X_t}\theta_t) = \mathbb{E}(\Delta_{t+1}|\sigma(\Delta_1, \dots, \Delta_t, \theta_0, \dots, \theta_t)), \tag{2.8}$$

so both conditional expectations coincide and as $\tau^{X_t}\theta_t = \omega_t$, they are functions of $\omega_t = (\boldsymbol{\omega}_t)_0$.

REMARK 2.6. Condition (A5) is slightly weaker than the corresponding conditions in [3, 22]. It is indeed equivalent to requiring that the set of points $z \in \Lambda$ with $\gamma_z > 0$ is not contained in any affine hyperplane of $\mathbb{R}$.

LEMMA 2.7. *Under Assumptions (A0), (A1) and (A5), the stationary Markov process $(\boldsymbol{\omega}_t)$ with distribution $\mathbf{P}$ is ergodic.*

PROOF. It suffices to prove that $\Pi h = h$ implies $h$ is $\mathbb{P}^e_\mu$-a.e. constant for each indicator function $h$. Consider a measurable set $C \subset \Omega$ such that $\Pi\mathbb{1}_C = \mathbb{1}_C$. Then, for each $\boldsymbol{\omega} = (\theta_t)_{t\in\mathbb{N}}$,

$$\mathbb{1}_C(\boldsymbol{\omega}) = (\Pi\mathbb{1}_C)(\boldsymbol{\omega}) = \sum_{z\in\Lambda} \pi_z(\theta_0)\mathbb{1}_C(\boldsymbol{\tau}^{1,z}\boldsymbol{\omega}),$$

which means that for each $z \in \Lambda$ such that $\gamma_z \neq 0$, $(\boldsymbol{\tau}^{1,z})^{-1}C \subset C$ holds $\mathbb{P}^e_\mu$-a.e. Then, by Assumption (A1), $\boldsymbol{\tau}^{1,z}C = C = (\boldsymbol{\tau}^{1,z})^{-1}C$ $\mathbb{P}^e_{\mu_e}$-a.e. since $\mathbb{P}^e_{\mu_e}$ is invariant under space–time translations. In addition, we will see shortly that the ellipticity Assumption (A5) implies that there exists $s \in \mathbb{N} \setminus \{0\}$ such that

$$\boldsymbol{\tau}^{s,0}C = C, \qquad \mathbb{P}^e_{\mu_e}\text{-a.s.} \tag{2.9}$$

The lemma thus follows by the ergodicity of $\mathbb{P}^e_{\mu_e}$ with respect to time translations which are multiples of $s$, which, in turn, follows from the mixing of the associated Markov process stated in Assumption (A0).

The proof of (2.9) is obvious if $\gamma_0 \neq 0$. To study the case $\gamma_0 = 0$, first notice that Assumption (A5) implies that the vectors in the set $V := \{(1,z) \in$ $\mathbb{R}\times\Lambda:\gamma_z\neq 0\}\subset\mathbb{R}^{d+1}$ must span $\mathbb{R}^{d+1}$. Otherwise, there would exist a vector $(a,l)\in\mathbb{Z}\times\mathbb{Z}^d$ such that $\langle(a,l),(1,z)\rangle=0$ for each $z\in V_0=\{(z\in\Lambda:\gamma_z\neq 0\}\subset\mathbb{R}^d$. But this means $\langle l,z\rangle=-a$ for each $z\in V_0$, which would contradict Assumption (A5).

Next, let $\{\bar{z}_i\}_{i=1}^{d+1}\subset V$ be a basis of $\mathbb{R}^{d+1}$. Accordingly, for each $s\in\mathbb{Z}$, we can solve the equation $\sum_{i=1}^{d+1}\alpha_i\bar{z}_i=s(1,0)$. In turn, this can be written in terms of a $(d+1)\times(d+1)$ invertible matrix with integer coefficients and a vector $\overline{\alpha}\in\mathbb{R}^{d+1}$ as $Z\overline{\alpha}=s(1,0)$. We choose $s\neq 0$ such that $\overline{\alpha}=sZ^{-1}(1,0)\in\mathbb{Z}^{d+1}$. If we let $A_+=\{i:\alpha_i>0\}$, then $\boldsymbol{\tau}^{\sum_{j=A_+}\alpha_j(1,z_j)}C=\boldsymbol{\tau}^{-\sum_{j\notin A_+}\alpha_j(1,z_j)}C$ $\mathbb{P}^e_{\mu_e}$-a.s., which implies (2.9) since for each $a\geq b$ and each set $A$, $(\boldsymbol{\tau}^{b,\zeta})^{-1}\boldsymbol{\tau}^{a,\eta}A\supset\boldsymbol{\tau}^{a-b,\eta-\zeta}A$. $\square$

Our first main result is a quenched CLT, that is, a CLT under the law $\mathbf{P}_{(\theta_t)}$, the measure $\mathbf{P}$ conditioned on the history $(\theta_t)_{t\in\mathbb{N}}$ of the environment.

THEOREM 1. *Under Assumptions (A0)–(A5), there exists a vector $v\in\mathbb{R}^d$ and a $d\times d$ matrix $\Sigma^2>0$ such that for $\mathbb{P}^e_{\mu_e}$-a.e. environment history $(\theta_t)\in\Omega$,*

$$\lim_{N\to\infty}\frac{1}{N}X_N=v,\qquad \mathbf{P}_{(\theta_t)}\text{-almost surely}\tag{2.10}$$

*and, letting $\widehat{X}_N:=X_N-Nv$,*

$$\frac{\widehat{X}_N}{\sqrt{N}}\Rightarrow\mathcal{N}(0,\Sigma^2)\qquad \text{under } \mathbf{P}_{(\theta_t)}.\tag{2.11}$$

PROOF. Recall Remarks 2.3 and 2.5 which allow us to interpret the random walk $(X_t)_{t\in\mathbb{N}}$ and all other random variables of interest as being defined on the probability space $(\boldsymbol{\Omega},\mathbf{P})$, and denote by $\mathbb{E}$ expectations with respect to $\mathbf{P}$. Consider the filtration

$$\mathcal{F}_t^0:=\sigma\{\Delta_1,\dots,\Delta_t,\theta_0,\dots,\theta_t\}.$$

Let $g(\omega_t)=\mathbb{E}(\Delta_{t+1}|\mathcal{F}_t^0)$. It is easy to see that this is a continuous local function of the process of the environment as seen from the particle (recall Section 2.2). We write

$$\begin{aligned}\Delta_{t+1}&=\Delta_{t+1}-\mathbb{E}(\Delta_{t+1}\mid\mathcal{F}_t^0)+\mathbb{E}(\Delta_{t+1}\mid\mathcal{F}_t^0)\\&=[\Delta_{t+1}-\mathbb{E}(\Delta_{t+1}\mid\mathcal{F}_t^0)]+g(\omega_t).\end{aligned}\tag{2.12}$$

Note that the first term is a martingale. Setting $v:=\mathbb{E}_\mu(g)$ and $g_0:=g-v$, we define

$$\widehat{\Delta}_{t+1}:=\Delta_{t+1}-v=\Delta_{t+1}-\mathbb{E}(\Delta_{t+1}\mid\mathcal{F}_t^0)+g_0(\omega_t).\tag{2.13}$$

Then, $\mathbb{E}(\widehat{\Delta}_{t+1}) = \mathbb{E}_\mu(g_0) = 0$ and as $\omega_t = (\boldsymbol{\omega}_t)_0$, $\sum_{t=0}^{N-1} \widehat{\Delta}_{t+1}$ is the sum of a martingale and an additive functional of the stationary ergodic process $(\boldsymbol{\omega}_t)$ under the law $\mathbf{P}$ (see Lemma 2.7). Hence, (2.10) follows from

$$\lim_{N\to\infty} N^{-1}(X_N - Nv) = \lim_{N\to\infty} \frac{1}{N} \sum_{t=0}^{N-1} (\Delta_{t+1} - \mathbb{E}(\Delta_{t+1} \mid \mathcal{F}_t^0) + g_0((\boldsymbol{\omega}_t)_0)) \\ = 0, \qquad \mathbf{P}\text{-a.s.}$$

Observe that "$\mathbf{P}_{(\theta_t)}$-almost surely for $\mathbb{P}_\mu^e$-a.e. environment history $(\theta_t) \in \Omega$" is the same as "$\mathbf{P}$-almost surely."

Of course, there is no such simple Fubini-type argument to pass from an unconditioned CLT (also known as an *annealed* CLT) to a conditional CLT. Nevertheless, we will first prove the unconditioned CLT since its proof is closely linked with a useful exponential estimate.

We wish to solve the equation $h - Sh = g_0$ that, thanks to Assumption (A2), has the bounded solution $h = \sum_{k=0}^\infty S^n g_0$. We can thus write (observing that $X_0 = 0$)

$$\begin{aligned} \widehat{X}_N &= \sum_{t=0}^{N-1} \widehat{\Delta}_{t+1} \\ &= \sum_{t=0}^{N-1} \{\Delta_{t+1} - \mathbb{E}(\Delta_{t+1} \mid \mathcal{F}_t^0) + h(\omega_{t+1}) - Sh(\omega_t)\} \\ &\quad + h(\omega_0) - h(\omega_N) \\ &= \sum_{t=0}^{N-1} \overline{\Delta}_{t+1} + h(\omega_0) - h(\omega_N), \end{aligned} \tag{2.14}$$

where $\overline{\Delta}_t := \Delta_t - \mathbb{E}(\Delta_t \mid \mathcal{F}_{t-1}^0) = h(\omega_t) - Sh(\omega_{t-1})$. If we let $M_n := \sum_{t=1}^n \overline{\Delta}_t$, then $M_n$ is an $\mathcal{F}_n^0$-martingale. Moreover, the $\overline{\Delta}_t$ are uniformly bounded random variables and they are almost surely functions of $\omega_t$ and $\omega_{t+1}$, so they are functions of $\boldsymbol{\omega}_t$ (see Remark 2.5). Therefore, $(\overline{\Delta}_{t+1}\overline{\Delta}_{t+1}^{\mathrm{T}})_t$ is a stationary and ergodic process and $N^{-1}\sum_{t=0}^{N-1} \overline{\Delta}_{t+1}\overline{\Delta}_{t+1}^{\mathrm{T}}$ converges almost surely to a symmetric matrix $\Sigma^2 \geq 0$. We note that this immediately implies the usual CLT

$$\frac{X_N - Nv}{\sqrt{N}} \Rightarrow \mathcal{N}(0, \Sigma^2) \qquad \text{under } \mathbf{P} \tag{2.15}$$

(see, e.g., [13], Theorem 3.2).

In addition, by a variant of Hoeffding's inequality for martingales (see, e.g., [12]), for sufficiently small $\varepsilon > 0$ and $L \in [0, N]$, the following holds:

$$\mathbf{P}\left\{\left|\frac{\widehat{X}_{N+L}}{\sqrt{(N+L)}} - \frac{\widehat{X}_N}{\sqrt{N}}\right| \geq \varepsilon\right\} \leq Ce^{-C\varepsilon^2(N/L)}. \tag{2.16}$$

This is easily seen, as follows. For sufficiently large $N$, we have, for each of the $d$ components separately,

$$\mathbf{P}\left\{\left|\frac{\widehat{X}_{N+L}}{\sqrt{(N+L)}}-\frac{\widehat{X}_N}{\sqrt{N}}\right|\geq\varepsilon\right\}\leq\mathbf{P}\left\{\frac{|M_{N+L}-M_N|}{\sqrt{(N+L)}}\geq\frac{\varepsilon}{4}\right\}$$
$$+\mathbf{P}\left\{|M_N|\left|\frac{1}{\sqrt{(N+L)}}-\frac{1}{\sqrt{N}}\right|\geq\frac{\varepsilon}{4}\right\}$$

and both terms can be estimated using Hoeffding's inequality, the first by $\exp(-C\varepsilon^2\frac{N}{L})$ and the second by $\exp(-C\varepsilon^2\frac{N^2}{L^2})$.

Next, we check that $\Sigma^2>0$. Indeed, if there exists $w\in\mathbb{R}^d$, $\|w\|=1$, such that $\langle w,\Sigma^2 w\rangle=0$, then

$$0=\mathbb{E}(\langle w,\Delta_{t+1}-\mathbb{E}(\Delta_{t+1}\mid\mathcal{F}_0^t)+h(\omega_{t+1})-Sh(\omega_t)\rangle^2),$$

which implies, for each $t\in\mathbb{N}$, $\langle w,\widehat{\Delta}_{t+1}\rangle=\langle w,h(\omega_t)-h(\omega_{t+1})\rangle$, hence

$$-\langle w,h(\omega_N)\rangle=\sum_{t=0}^{N-1}\langle w,\widehat{\Delta}_{t+1}\rangle+\langle w,h(\omega_0)\rangle. \tag{2.17}$$

This is in contradiction with the boundedness of $h$. In fact, on the one hand, (2.17) implies $|\sum_{t=0}^{N-1}\langle w,\widehat{\Delta}_{t+1}\rangle|\leq 2\|h\|_\infty$. On the other hand, by Assumption (A5), there exists a probability larger than $\gamma^N$ to have $|\sum_{t=0}^{N-1}\langle w,\Delta_{t+1}-v\rangle|\geq C_w N$.

To obtain more refined information, it is convenient to consider the finer filtration

$$\mathcal{F}_t:=\sigma\{\Delta_1,\ldots,\Delta_t,(\theta_s)_{s\in\mathbb{N}}\}$$

and the decomposition

$$\widehat{\Delta}_{t+1}=(\widehat{\Delta}_{t+1}-\mathbb{E}(\widehat{\Delta}_{t+1}\mid\mathcal{F}_t))+\mathbb{E}(\widehat{\Delta}_{t+1}\mid\mathcal{F}_t). \tag{2.18}$$

Clearly, $Z_t:=\sum_{s=1}^t\widehat{\Delta}_s-\mathbb{E}(\widehat{\Delta}_s\mid\mathcal{F}_{s-1})$ is a martingale with respect to the filtration $\mathcal{F}_t$. Let

$$\tilde{g}(\theta):=\sum_{z\in\Lambda}z\pi_z(\theta_0)-v,$$

so $\tilde{g}$ is a continuous local function on $\Theta$. Then,

$$\mathbb{E}(\widehat{\Delta}_{t+1}\mid\mathcal{F}_t)=\sum_{z\in\Lambda}z\pi_z(\tau^{X_t}\theta_t)-v=\tilde{g}(\tau^{X_t}\theta_t)=\tilde{g}(\omega_t). \tag{2.19}$$

Recall that $\boldsymbol{\omega}=(\omega_s)_{s\in\mathbb{N}}$. Define $G(\boldsymbol{\omega}):=\tilde{g}(\omega_0)$. By Lemma 2.1, Remark 2.3 and Section 2.2 on the environment as seen from the particle, we have

$$\begin{aligned}\mathbb{E}(\widehat{\Delta}_{t+1}\mid\mathcal{F}_t)&=\tilde{g}(\omega_t)=G(\boldsymbol{\omega}_t),\\ \mathbb{E}(\widehat{\Delta}_{t+1}\mid\mathcal{F}_0)&=\Pi^t G(\boldsymbol{\omega}_0)=\mathbb{E}(G(\boldsymbol{\omega}_t)|\sigma(\boldsymbol{\omega}_0)),\end{aligned} \tag{2.20}$$

where $(\boldsymbol{\omega}_t)$ is the Markov process defined in Section 2.2. Thus, the remainder in (2.18) is an additive functional of this Markov process.

The next idea, following [20], is to use [19], Theorem 1, to conclude. To be precise, [20] uses [19] in conjunction with the theory of fractional coboundaries developed in [8]. In fact, since we are discussing random walks with bounded increments, the use of [8] is not really necessary and a slightly more quantitative version of [19] allows us to conclude by a simple Borel–Cantelli argument; see [21] for a similar strategy. For the reader's convenience we present the needed modifications of the arguments from [19, 21] in the Appendix. Indeed, Theorem A.2 shows the following. Given any $w \in \mathbb{R}^d$ and $\rho \geq 0$, if

$$\sum_{n=1}^{\infty} n^{-3/2}(\ln n)^{\rho} \left\| \sum_{k=0}^{n-1} \langle w, \Pi^k G \rangle \right\|_{L^2(\Omega, \mathbb{P}^e_\mu)} < \infty, \tag{2.21}$$

then $\sum_{s=0}^{t} G(\boldsymbol{\omega}_s)$ can be decomposed, under the stationary law $\mathbf{P}$ of $(\boldsymbol{\omega}_t)_{t\in\mathbb{N}}$, as $\tilde{M}_t + R_t$, where $\tilde{M}_t$ is an $L^2(\boldsymbol{\Omega}, \mathbf{P})$ martingale with respect to the filtration $\mathcal{F}_t = \sigma\{\boldsymbol{\omega}_0, \ldots, \boldsymbol{\omega}_t\}$ (see Remark 2.3) and $\lim_{N\to\infty} N^{-1/2} R_N = 0$ in $L^2(\boldsymbol{\Omega}, \mathbf{P})$. In addition, $\tilde{M}_t - \tilde{M}_{t-1}$ is $\sigma\{\boldsymbol{\omega}_{t-1}, \boldsymbol{\omega}_t\}$-measurable and it can be written as $H(\boldsymbol{\omega}_{t-1}, \boldsymbol{\omega}_t)$ for some $\mathbb{R}^d$-valued function $H \in L^2(\Theta^2, \mathbf{P}_2)$, where $\mathbf{P}_2$ is the two dimensional marginal of $\mathbf{P}$. That is, $\int_{\Theta^2} f(\omega, \omega')\mathbf{P}_2(d\omega, d\omega') := \sum_{z\in\Lambda} \int_{\Theta^2} \pi_z(\omega) f(\omega, \boldsymbol{\tau}^{1,z}\omega)\mathbb{P}^e_\mu(d\omega)$. We thus have that

$$\hat{X}_t = \sum_{s=1}^{t} \hat{\Delta}_s = Z_t + \tilde{M}_{t-1} + R_{t-1} = Z_t + \tilde{M}_t + \tilde{R}_t,$$

where $Z_t + \tilde{M}_t$ is an $\mathcal{F}_t$-martingale and $\tilde{R}_t = R_{t-1} - H(\boldsymbol{\omega}_{t-1}, \boldsymbol{\omega}_t)$ is of order $\mathcal{O}(t^{1/2}(\ln t)^{-\rho})$ in $L^2(\boldsymbol{\Omega}, \mathbf{P})$; see Theorem A.2. Define the $\mathbb{R}^{d\times d}$-valued function $F \in L^1(\Omega, \mathbb{P}^e_\mu)$ by

$$F(\omega_s) := \mathbb{E}((Z_{s+1} + \tilde{M}_{s+1} - Z_s - \tilde{M}_s)(Z_{s+1} + \tilde{M}_{s+1} - Z_s - \tilde{M}_s)^{\mathrm{T}} | \mathcal{F}_s)$$

(observe that $F$ depends only on $\omega_s$, due to Remark 2.3). Then, the average (conditional) quadratic variation of the $\mathbb{R}^d$-valued martingale $Z_t + \tilde{M}_t$, for $\mathbb{P}^e_\mu$-a.e. environment history $(\theta_t)_{t\in\mathbb{N}}$, is given by

$$\begin{aligned} &\lim_{t\to\infty} \frac{1}{t} \sum_{s=0}^{t-1} \mathbb{E}((Z_{s+1} + \tilde{M}_{s+1} - Z_s - \tilde{M}_s)(Z_{s+1} + \tilde{M}_{s+1} - Z_s - \tilde{M}_s)^{\mathrm{T}} | \mathcal{F}_s) \\ &\quad = \lim_{t\to\infty} \frac{1}{t} \sum_{s=0}^{t-1} F(\boldsymbol{\omega}_s) = \mathbb{E}_{\mathbb{P}^e_\mu}(F), \end{aligned}$$

where we have used Birkhoff's theorem and the ergodicity of the process $(\boldsymbol{\omega}_t)$ under $\mathbf{P}$ (see Lemma 2.7). Hence, for $\mathbb{P}^e_\mu$-a.e. $(\theta_t)$, we have convergence $t^{-1/2}(Z_t + \tilde{M}_t) \Rightarrow \mathcal{N}(0, \mathbb{E}_{\mathbb{P}^e_\mu}(F))$, by standard martingale CLT convergence theorems.

Indeed, one may apply [13], Theorem 3.2, to the conditional martingales "$Z_t + \tilde{M}_t$ given $(\theta_t)$." To do this, one needs to check that for $\mathbb{P}^e_\mu$-a.e. $(\theta_t)$, this conditional martingale satisfies a conditional Lindeberg condition which, in turn, is implied by the slightly stronger requirement that $\lim_{t\to\infty} \frac{1}{t}\sum_{s=0}^{t-1} \mathbb{E}(f_{s\varepsilon^2}(\boldsymbol{\omega}_s, \boldsymbol{\omega}_{s+1})|\mathcal{F}_s) = 0$ $\mathbf{P}$-almost surely for each $\varepsilon > 0$, where $f_u(\boldsymbol{\omega}_s, \boldsymbol{\omega}_{s+1}) = \langle \xi_{s+1}, \xi_{s+1}\rangle \mathbb{1}_{\{\langle \xi_{s+1}, \xi_{s+1}\rangle > u\}}$ and $\xi_{s+1} = Z_{s+1} + \tilde{M}_{s+1} - Z_s - \tilde{M}_s$. But, for each $u > 0$,

$$\begin{aligned}\lim_{t\to\infty} \frac{1}{t}\sum_{s=0}^{t-1} \mathbb{E}(f_{s\varepsilon^2}(\boldsymbol{\omega}_s, \boldsymbol{\omega}_{s+1})|\mathcal{F}_s) &\le \lim_{t\to\infty} \frac{1}{t}\sum_{s=0}^{t-1} \mathbb{E}(f_u(\boldsymbol{\omega}_s, \boldsymbol{\omega}_{s+1})|\mathcal{F}_s) \\ &= \mathbb{E}(f_u) = \mathbb{E}(\langle \xi_{s+1}, \xi_{s+1}\rangle \mathbb{1}_{\{\langle \xi_{s+1}, \xi_{s+1}\rangle > u\}})\end{aligned}$$

$\mathbf{P}$-almost surely by Birkhoff's theorem and the observation that $(\boldsymbol{\omega}_t)$ is a Markov process, and this value tends to zero as $u \to \infty$.

Note also that

$$\mathbb{E}_{\mathbb{P}^e_\mu}(F) = \lim_{N\to\infty} N^{-1}\mathbb{E}((Z_N + \tilde{M}_N)^2) = \lim_{N\to\infty} N^{-1}\mathbb{E}(\widehat{X}_N^2) = \Sigma^2, \tag{2.22}$$

in particular, $\mathbb{E}_{\mathbb{P}^e_\mu}(F) = \Sigma^2 > 0$.

The last task is to prove that the remainder $t^{-1/2}\tilde{R}_t$ converges to zero almost surely. Given the available estimates, we first prove it only for the subsequence $t \in T_a := \{[1 + jk^{-a}]2^k\}_{k\in\mathbb{N}, 0\le j<k^a}$, where $a > 1$ is such that $2\rho > 1 + a$. Here, we assume that $\rho > 1$, for which (2.21) holds. Indeed, by Theorem A.2 and Chebyshev's inequality, it follows that

$$\sum_{t\in T_a} \mathbf{P}(\{|t^{-1/2}\tilde{R}_t| \ge \varepsilon\}) \le C\varepsilon^{-2}\sum_{t\in T_a}(\ln t)^{-2\rho} \le C\varepsilon^{-2}\sum_{k\in\mathbb{N}} k^{a-2\rho} < \infty.$$

By Borel–Cantelli, it follows that $t^{-1/2}\tilde{R}_t$ converges to zero almost surely along the subsequence $T_a$. Accordingly, along the subsequence $T_a$, conditioned on $\mathbb{P}^e_\mu$-almost every environment history, the random variables $N^{-1/2}\widehat{X}_N$ converge weakly to a Gaussian with variance $\Sigma^2$. To conclude, we use the fact, quantified in (2.16), that $N^{-1/2}\widehat{X}_N$ changes very slowly. Given any $n \in \mathbb{N}$, let $n_{T_a} \in T_a$ be the element of $T_a$ closest to $n$. Then, (2.16) implies

$$\sum_{n\in\mathbb{N}} \mathbf{P}\left(\left\{\left|\frac{\widehat{X}_n}{\sqrt{n}} - \frac{\widehat{X}_{n_{T_a}}}{\sqrt{n_{T_a}}}\right| \ge \varepsilon\right\}\right) \le \sum_{n\in\mathbb{N}} Ce^{-C\varepsilon^2 n/|n-n_{T_a}|} \le \sum_{k\in\mathbb{N}} C2^k e^{-C\varepsilon^2 k^a} < \infty,$$

provided $a > 1$. Hence, again by Borel–Cantelli, the sequence $(\frac{\widehat{X}_n}{\sqrt{n}} - \frac{\widehat{X}_{n_{T_a}}}{\sqrt{n_{T_a}}})$ converges to zero $\mathbb{P}^e_\mu$-almost surely, which implies the claimed result.

The theorem is thus proved, provided (2.21) holds with $\rho > 1$. In Section 3.2, we will see, following [5], that such an estimate is equivalent to estimating the number of times two independent walks in the same environment come close. We will then show that (2.21) is indeed satisfied under Assumptions (A0)–(A5). □

2.4. *A concrete model*: *weakly coupled Markov chains.* We have seen that under some assumptions, it is possible to prove a quenched CLT theorem for the random walk. It is now time to present a concrete class of examples in which such assumptions are satisfied.

Let $K(\theta, dy)$ be a transition kernel that specifies the transition probability from $\theta^0 \in I$ to $y \in I$ given the rest of the configuration $\theta^{\neq 0} := (\theta^p)_{p \neq 0}$. Clearly, $\int_I K(\theta, dy) = 1$. We further require that for each $u \in \mathcal{C}^0(I)$ and $q \in \mathbb{Z}^d$, the function $\tilde{u}$ defined by $\tilde{u}(\theta) := \int_I u(y) K(\theta, dy)$ belongs to $\mathcal{C}^0(\Theta)$. So, we can define a Feller Markov operator $\mathcal{K} : \mathcal{C}^0(\Theta) \to \mathcal{C}^0(\Theta)$,

$$(\mathcal{K} f)(\theta) := \int_\Theta \prod_{q \in \mathbb{Z}^d} K(\tau^q \theta, dy_q) f(y). \tag{2.23}$$

In fact, $\mathcal{K}$ is clearly well defined on $\mathcal{C}^0_{\mathrm{loc}}(\Theta)$ and it extends, by continuity, to all of $\mathcal{C}^0(\Theta)$.

Assumption (A6) (*Local mixing*). For each $\theta, \tilde{\theta} \in \Theta$ such that $\theta^q = \tilde{\theta}^q$ for all $q \neq 0$, we assume

$$|K(\theta, dy) - K(\tilde{\theta}, dy)| \leq 2 d_0,$$

where the norm refers to the total variation of measures.

Assumption (A7) (*Weak coupling*). For each $p \neq 0$ and $\theta, \tilde{\theta} \in \Theta$ such that $\theta^q = \tilde{\theta}^q$ for each $q \neq p$, we assume

$$|K(\theta, dy) - K(\tilde{\theta}, dy)| \leq d_q.$$

Assumption (A8) (*Long range bound*). There exists $\xi' > \max\{4 + d, 2d\}$ such that for all $L > 0$,

$$\sum_{\|q\| \geq L} d_q \leq C L^{-\xi'}.$$

Assumption (A9) (*Dobrushin-like conditions*). Write the transition probabilities of the random walk in the form $\pi_z(\theta) = a_z + \hat{\pi}_z(\theta)$, where $a_z \in [0,1]$, $\sum_{z \in \Lambda} a_z = 1$ and $\hat{\pi}_z$ depends only on $\theta^p$, $p \in \Lambda$. Then, setting $D := \sum_{z \in \Lambda} \|\hat{\pi}_z\|_\infty$, consider the following hierarchy of conditions:

(a) $\eta_0 := \sum_{q \in \mathbb{Z}^d} d_q < 1$;
(b) $\eta_1 := (1 + (1 + 2|\Lambda|) D) \eta_0 < 1$;
(c) $\eta_1 < 1$; $D < 1$; $\eta_0 < (1 - D)^{(\xi'(1+d))/(\xi' - d)}$.

REMARK 2.8. In the case of finite-range interactions ($d_q = 0$ for $\|q\|$ larger than some $R > 0$), a polynomial decay of time correlations suffices to prove the CLT. One can then use the strategy employed at the end of Section 4.5 to obtain the result under a weaker smallness condition than the one stated in Assumption (A9).

In Section 4, we will prove the following theorem.

THEOREM 2. (i) *Each Markov environment satisfying Assumptions (A6), (A7) and* (A9)(a) *enjoys property (A0), although only relative to* $\mu_e$.

(ii) *If* (A9)(b) *is also satisfied, then (A0) also holds for* $\mu$.

(iii) *If, in addition, it satisfies property (A8), then it also satisfies Assumptions (A2), (A3) and (A4).*

(iv) *Finally, if Assumptions* (A9)(c) *is also satisfied, then Assumptions (A1) holds.*

REMARK 2.9. The above assumptions are very similar to (although much more general than the ones used in [1, 3]. In fact, in [1, 3], the Markov chains are independent at each site, this corresponding to $d_q = 0$ for each $q \neq 0$ so that $\eta_0 = d_0$ in Assumption (A9). Also, the random walk is a nearest neighbor walk and the transition probabilities depend only on the site presently visited. This corresponds to having $\pi_z$ depending only on $\theta^0$ in the present setting, which would imply that the constant $(1 + 2|\Lambda|)$ in Assumption (A9) can be replaced by 3. In the situation described above, one can compare our results with the previous ones (keeping in mind that [3] holds only for $d \geq 3$, while [1] only for $d > 7$). For example, let us compare with [1], which presents sharper results in its realm of applicability.

The assumptions of [1] are in terms of the two parameters $\kappa$ and $\varepsilon$, where $1 - \kappa$ gives a bound for the rate of mixing and thus corresponds to our $\eta_0$. The parameter $\varepsilon$ is the best constant such that $\pi_z \geq \varepsilon a_z$. The equivalent of Assumption (A9) (the only relevant condition here) in [1] reads $\kappa + \varepsilon^2 > 1$, that is, $\kappa = 1 - \varepsilon^\gamma$ for some $\gamma > 2$. This is sufficient for the annealed CLT in [1], but $\gamma > 6$ is assumed for the quenched CLT (indeed, not only CLTs but Donsker-type invariance principles are proved in [1]). To see the relation with our conditions, let us consider the simple case in which $|\hat{\pi}_z| \leq (1 - \varepsilon) a_z$. [In fact, the following holds more generally, as can be seen by using the decomposition (4.15).] Then, $D = 1 - \varepsilon$ and Assumption (A9)(b), in the language of [1], reads $(1 + 3(1 - \varepsilon))(1 - \kappa) < 1$, that is, $\kappa + (4 - 3\varepsilon)^{-1} > 1$. One can easily verify that the above condition is better than $\kappa + \varepsilon^2 > 1$, provided $\varepsilon \leq 0.75$. Furthermore, if $\kappa = 1 - \varepsilon^\gamma$ for some $\gamma \geq 3$, then $\kappa + (4 - 3\varepsilon)^{-1} > 1$ for all $0 \leq \varepsilon < 1$. As for the condition Assumption (A9)(c), a direct computation along the lines of Section 4.5 yields that it can be replaced, in this case, by $\eta_0 + D < 1$, that is, $\kappa + \varepsilon > 1$, which is always weaker than the above. The reason is that all truncated operators in the proof of Lemma 4.8 coincide in this case with the untruncated ones so that the estimates in (4.21) and (4.27) can be replaced by equalities. As a result, the estimate in Lemma 4.8 is uniform in $n$. In other words, our conditions are weaker than the condition under which the quenched invariance principle is proved in [1].

**3. Proofs: CLT under Assumptions (A0)–(A5).** In this section, we prove Theorem 1.

3.1. *Equivalence of* (2.21) *with a two-walks estimate.* As discussed in Section 2.3, it suffices to prove (2.21). Recalling (2.20), for each $w \in \mathbb{R}^d$, $\|w\| \leq 1$,

$$\left\| \sum_{t=0}^{N-1} \langle w, \Pi^t G \rangle \right\|_2^2 = \sum_{t,s=1}^{N} \mathbb{E}(\mathbb{E}(\langle w, \widehat{\Delta}_t \rangle | \mathcal{F}_0) \mathbb{E}(\langle w, \widehat{\Delta}_s \rangle \mid \mathcal{F}_0)).$$

The above formula has a very interesting interpretation: consider two independent random walks $X_n, Y_n$, both starting from zero and evolving in the same environment $(\theta_t)$ described by the transition probabilities (2.1). That is, setting $\widehat{\Delta}_{t+1}^X := X_{t+1} - X_t - v$, $\widehat{\Delta}_{t+1}^Y := Y_{t+1} - Y_t - v$, we have

$$\begin{aligned} P(\{(\widehat{\Delta}_{t+1}^X, \widehat{\Delta}_{t+1}^Y) = (z, z'), \theta_{t+1} \in A\} \mid X_t, Y_t, \theta_t) \\ = \pi_z(\tau^{X_t} \theta_t) \pi_{z'}(\tau^{Y_t} \theta_t) p(\theta_t, A). \end{aligned}$$

Let us call $\mathbb{P}_\nu^2$ the law of such a process when the environment is started with the measure $\nu$ and denote by $\mathbb{E}_\nu^2$ the corresponding expectation. Then,

$$\left\| \sum_{t=0}^{N-1} \langle w, \Pi^t G \rangle \right\|_2^2 = \sum_{t,s=1}^{N} \mathbb{E}_\mu^2(\langle w, \widehat{\Delta}_t^X \rangle \langle w, \widehat{\Delta}_s^Y \rangle). \tag{3.1}$$

Remark 3.1. Note that if the process $(X_n, Y_n)$ satisfies the CLT (which is, in fact, a consequence of what we will prove later on), then (3.1) corresponds to the off-diagonal part of the covariance of such a process. From this point of view, condition (2.21) says that the two walks are asymptotically independent.

3.2. *The two-walks estimate: off-diagonal variance.* Let $L_N := A \ln N$, for some fixed $A > 0$ to be chosen later [see (3.10) and (3.12)].

Lemma 3.2. *There exists $\delta > 0$ such that*

$$\mathbb{E}_\mu^2(\operatorname{Card}\{t \leq N : \|X_t - Y_t\| \leq L_N\}) \leq C N^{1-\delta} \qquad (N \in \mathbb{N}). \tag{3.2}$$

This lemma is proved in Section 3.3. We now complete the proof of (2.21) using the lemma.

Let us introduce the filtrations

$$\mathcal{F}_t^{0,XY} := \sigma\{X_0, \dots, X_t, Y_0, \dots, Y_t, \theta_0, \dots, \theta_t\},$$
$$\mathcal{F}_t^{XY} := \sigma\{X_0, \dots, X_t, Y_0, \dots, Y_t, (\theta_s)_{s\in\mathbb{N}}\}$$

and the filtrations $\mathcal{F}_t^{0,X}$ and $\mathcal{F}_t^{0,Y}$, which are just the '$X$- and $Y$-versions' of the previously introduced $\mathcal{F}_t^0$. To estimate (3.1), we start by considering the case $t < s$. As $\widehat{\Delta}_t^X$ and $\widehat{\Delta}_s^Y$ are conditionally independent given $\mathcal{F}_t^{0,XY}$, we have

$$\begin{aligned}\mathbb{E}^2_\mu(\langle w, \widehat{\Delta}_t^X\rangle\langle w, \widehat{\Delta}_s^Y\rangle) &= \mathbb{E}^2_\mu(\mathbb{E}^2_\mu(\langle w, \widehat{\Delta}_t^X\rangle \mid \mathcal{F}_t^{0,XY})\mathbb{E}^2_\mu(\langle w, \widehat{\Delta}_s^Y\rangle \mid \mathcal{F}_t^{0,XY}))\\ &= \mathbb{E}^2_\mu(\langle w, \mathbb{E}^2_\mu(\widehat{\Delta}_t^X \mid \mathcal{F}_t^{0,X})\rangle\langle w, \mathbb{E}^2_\mu(\widehat{\Delta}_s^Y \mid \mathcal{F}_t^{0,Y})\rangle).\end{aligned}$$

Calling $\omega^Y$ the environment as seen from $Y$, we have

$$\mathbb{E}^2_\mu(\widehat{\Delta}_s^Y \mid \mathcal{F}_t^{0,Y}) = S^{s-t-1}g(\omega_t^Y).$$

Assumption (A2) then implies

$$|\mathbb{E}^2_\mu(\langle w, \widehat{\Delta}_t^X\rangle\langle w, \widehat{\Delta}_s^Y\rangle)| \le C\eta^{s-t-1}.$$

Hence, for each $a > \frac{\xi}{2}$, setting $b = \frac{2(\rho+a)}{\ln\eta^{-1}}$ and defining $T_N := b\ln(\ln N)$, we have

$$\begin{aligned}\sum_{\substack{|t-s|\ge T_N\\ 1\le t,s\le N}} |\mathbb{E}^2_\mu(\langle w, \widehat{\Delta}_t^X\rangle\langle w, \widehat{\Delta}_s^Y\rangle)| &\le C_a N(\ln N)^{-2(\rho+a)},\\ &\le C_a N(\ln N)^{-2\rho-\xi}.\end{aligned} \tag{3.3}$$

Since the roles of $t$ and $s$ are interchangeable, it remains to consider the cases for which $s \ge t \ge T_N$ and $s - t \le T_N$. Let us write $\widehat{\Delta}_{w,t}^X = \langle w, \widehat{\Delta}_t^X\rangle$ and $\widehat{\Delta}_{w,t}^Y = \langle w, \widehat{\Delta}_t^Y\rangle$. We can then write

$$\begin{aligned}|\mathbb{E}^2_\mu(\langle w, \widehat{\Delta}_t^X\rangle\langle w, \widehat{\Delta}_s^Y\rangle)| &\le |\mathbb{E}^2_\mu(\mathbb{1}_{\{\|X_{t-T_N}-Y_{t-T_N}\|\ge L_N\}}\widehat{\Delta}_{w,t}^X\widehat{\Delta}_{w,s}^Y)|\\ &\quad + C\mathbb{P}^2_\mu(\{\|X_{t-T_N} - Y_{t-T_N}\| < L_N\}).\end{aligned} \tag{3.4}$$

Let us set $A_{N,t} := \{\|X_{t-T_N} - Y_{t-T_N}\| \ge L_N\}$. Using the fact that $\widehat{\Delta}_{w,t}^X$ and $\widehat{\Delta}_{w,s}^Y$ are conditionally independent given $\mathcal{F}_{t-T_N}^{XY}$, we can write

$$\begin{aligned}&|\mathbb{E}^2_\mu(\mathbb{1}_{A_{N,t}}\mathbb{E}^2_\mu(\widehat{\Delta}_{w,t}^X\widehat{\Delta}_{w,s}^Y \mid \mathcal{F}_{t-T_N}^{0,XY}))|\\ &\quad = |\mathbb{E}^2_\mu(\mathbb{1}_{A_{N,t}}\mathbb{E}^2_\mu[\mathbb{E}^2_\mu(\widehat{\Delta}_{w,t}^X \mid \mathcal{F}_{t-T_N}^{XY})\mathbb{E}^2_\mu(\widehat{\Delta}_{w,s}^Y \mid \mathcal{F}_{t-T_N}^{XY}) \mid \mathcal{F}_{t-T_N}^{0,XY}])|.\end{aligned}$$

We want to estimate the conditional expectation with respect to $\mathcal{F}^{0,XY}_{t-T_N}$. To this end, we fix $\theta_0,\dots,\theta_{t-T_N}$, $X_0,\dots,X_{t-T_N}$ and $Y_0,\dots,Y_{t-T_N}$. Then, $\mathbb{E}^2_\mu(\widehat{\Delta}^X_{w,t}\mid\mathcal{F}^{XY}_{t-T_N})$ and $\mathbb{E}^2_\mu(\widehat{\Delta}^Y_{w,s}\mid\mathcal{F}^{XY}_{t-T_N})$ are functions of two subsets of the spatial coordinates of the $T_N+s-t\le 2T_N$ variables $\theta_{t-T_N+1},\dots,\theta_s$ and these subsets are separated by a distance $L_N-CT_N\ge L_N/2$. Since $\mathbb{E}^2_\mu(f(\theta_{t-T_N},\dots,\theta_s)\mid\mathcal{F}^{0,XY}_{t-T_N})=\mathbb{E}_{\mathbb{P}^e_{\theta_{t-T_N}}}(f)$, we can apply Assumption (A4) and estimate the first term on the right-hand side of (3.4) by

$$
\begin{aligned}
&|\mathbb{E}^2_\mu(\mathbb{1}_{A_{N,t}}\mathbb{E}^2_\mu(\widehat{\Delta}^X_{w,t}\widehat{\Delta}^Y_{w,s}\mid\mathcal{F}^{0,XY}_{t-T_N}))|\\
&\qquad\le|\mathbb{E}^2_\mu(\mathbb{1}_{A_{N,t}}\mathbb{E}^2_\mu(\widehat{\Delta}^X_{w,t}\mid\mathcal{F}^{0,XY}_{t-T_N})\mathbb{E}^2_\mu(\widehat{\Delta}^Y_{w,s}\mid\mathcal{F}^{0,XY}_{t-T_N}))|+CL_N^{-\xi}T_N^{\tilde{\xi}}\\
&\qquad\le\mathbb{E}(|S^{T_N}g(\omega^X_{t-T_N})|\cdot|S^{s-t+T_N}g(\omega^Y_{t-T_N})|)+CL_N^{-\xi}T_N^{\tilde{\xi}}\\
&\qquad\le C(\eta^{T_N}+L_N^{-\xi}T_N^{\tilde{\xi}}).
\end{aligned}
\tag{3.5}
$$

Hence, if $\xi>2\rho$, then

$$
\begin{aligned}
&\sum_{\substack{|t-s|\le b\ln(\ln N)\\ 1\le t,s\le N}}|\mathbb{E}^2_\mu(\langle w,\widehat{\Delta}^X_t\rangle\langle w,\widehat{\Delta}^Y_l\rangle)|\\
&\qquad\le\frac{C_{A,b}N(\ln\ln N)^{\tilde{\xi}+1}}{(\ln N)^{\xi}}+T_N\sum_{t=T_N}^{N-1}\mathbb{P}^2_\mu(A^c_{N,t})+T_N^2.
\end{aligned}
\tag{3.6}
$$

Combining (3.3) and (3.6), proves (2.21), provided that $\xi>2\rho+2$. [Note that because $\xi>4$ in Assumption (A4), we may choose $\rho>1$, as required at the end of Section 2.3.]

3.3. *Estimating the number of close encounters.* We first reduce (3.2) to a simpler inequality.

Lemma 3.3. *There exist $\beta\in(0,1),C_0>0$ such that for any $\theta\in\Theta$ and any $a,b$ such that $\|a-b\|>L_N$, we have*

$$
\mathbb{P}^2_\theta(\|X_j-Y_j\|>L_N\ \text{for}\ j=1,2,\dots,N|X_0=a,Y_0=b)\ge\frac{C_0}{N^\beta}.
\tag{3.7}
$$

*[Here, $\mathbb{P}^2_\theta$ is the underlying probability for the process $(\theta_t,X_t,Y_t)$ started at $\theta_0=\theta$.]*

Proof of Lemma 3.2. We start by noticing that Assumption (A5) implies, for each $a,b\in\mathbb{Z}^d$, $\|a-b\|\le L_N$, that

$$
\mathbb{P}^2_\theta\biggl(\biggl\{\sup_{0\le i\le L_N}\|X_i-Y_i\|\ge L_N\biggr\}\biggm| X_0=a,Y_0=b\biggr)\ge\gamma^{L_N},
\tag{3.8}
$$

the latter being the probability of one fixed path in which $X_i, Y_i$ get further and further apart at each step. Accordingly, for each $\varrho < 1-\beta$, we have

$$
\begin{aligned}
&\mathbb{P}^2_\theta\left(\left\{\sup_{0\le i\le N^\varrho}\|X_i - Y_i\| \le L_N\right\}\Big| X_0 = a, Y_0 = b\right)\\
&\quad\le \prod_{j=1}^{N^\varrho L_N^{-1}} (1-\gamma^{L_N}) \qquad (3.9)\\
&\quad\le e^{-\gamma^{L_N} L_N^{-1} N^\varrho} \le e^{-N^{\varrho - 2A\ln\gamma^{-1}}} \le e^{-N^{\varrho/2}},
\end{aligned}
$$

where we have chosen $A$ such that

$$\varrho > 4A\ln\gamma^{-1}. \qquad (3.10)$$

Next, consider the sets $B_R^- := \{(x,y) : \|x-y\| \le R\}$, $B_R^+ := \{\|x-y\| > R\}$ and the stopping times, for $k > 0$,

$$
\begin{aligned}
s_0 &:= \inf\{j \in \mathbb{N} : (X_j, Y_j) \in B^-_{L_N}, (X_{j+1}, Y_{j+1}) \in B^+_{L_N}\},\\
s_{2k} &:= \inf\{j \in \mathbb{N} : j > s_{2k-2}, (X_j, Y_j) \in B^-_{L_N}, (X_{j+1}, Y_{j+1}) \in B^+_{L_N}\},\\
s_1 &:= \inf\{j \in \mathbb{N} : j > s_0, (X_j, Y_j) \in B^+_{L_N}, (X_{j+1}, Y_{j+1}) \in B^-_{L_N}\},\\
s_{2k+1} &:= \inf\{j \in \mathbb{N} : j > s_{2k-1}, (X_j, Y_j) \in B^+_{L_N}, (X_{j+1}, Y_{j+1}) \in B^-_{L_N}\}.
\end{aligned}
$$

Clearly, $s_{2k} < s_{2k+1} < s_{2k+2}$ and $s_k > k$. As $X_0 = Y_0$, these stopping times are adapted to the filtration $\mathcal{F}^{0,XY}_t$. With this notation, (3.9) implies

$$\mathbb{P}^2_\theta\left(\left\{\sup_{i\le N}(s_{2i} - s_{2i-1}) > N^\varrho\right\}\right) \le N \sup_{i\le N} \mathbb{P}^2_\theta(\{s_{2i} - s_{2i-1} > N^\varrho\}) \le N e^{-N^{\varrho/2}}.$$

Let us set $J := \inf\{k \in \mathbb{N} : s_{2k+1} \ge N\} + 1$. Obviously, $J \le N/2 + 1$. Clearly, $J$ is the number of intervals in which the two walks are closer than $L_N$ before time $N$. Since the above estimate tells us that such intervals are shorter than $N^\varrho$, with overwhelming probability, we have

$$\mathbb{E}^2_\theta(\mathrm{Card}\{n < N : \|X_n - Y_n\| \le L_N\}) \le N^2 e^{-N^{\varrho/2}} + N^\varrho \mathbb{E}^2_\theta(J). \qquad (3.11)$$

It remains to investigate the lengths of the intervals of time in which the two walks are closer than $L_N$ or (which is the same) the ones in which they are further apart than $L_N$. Let $S_n := \{\sup_{k\le n}(s_{2k+1} - s_{2k}) < N\}$ and denote by $\mathcal{F}^{0,XY}_{s_{2k}}$ the $\sigma$-algebra associated to the filtration $\mathcal{F}^{0,XY}_t$ and the stopping time $s_{2k}$. Then,

$$
\begin{aligned}
\mathbb{P}^2_\theta(\{J > n+1\}) &= \mathbb{P}^2_\theta(\{s_{2k+1} < N\ \forall k \le n\})\\
&\le \mathbb{P}^2_\theta(\{s_{2k+1} - s_{2k} < N\ \forall k \le n\}) = \mathbb{E}^2_\theta(\mathbb{1}_{S_n}).
\end{aligned}
$$

Thus, by (3.7),

$$\begin{aligned}\mathbb{P}^2_\theta(\{J > n+1\}) &\le \mathbb{E}^2_\theta(\mathbb{1}_{S_n}) = \mathbb{E}^2_\theta(\mathbb{1}_{S_{n-1}}\mathbb{P}^2_\theta(\{s_{2n+1} - s_{2n} < N\} \mid S_{n-1})) \\ &= \mathbb{E}^2_\theta(\mathbb{1}_{S_{n-1}}\mathbb{P}^2_\theta(\{s_{2n+1} - s_{2n} < N\} \mid \mathcal{F}^{0,XY}_{s_{2n}})) \\ &\le \left(1 - \frac{C_0}{N^\beta}\right)\mathbb{E}^2_\theta(\mathbb{1}_{S_{n-1}}) \le \cdots \le \left(1 - \frac{C_0}{N^\beta}\right)^n.\end{aligned}$$

Thus, letting $1 - \varrho > \alpha > \beta$, it follows that

$$\mathbb{P}^2_\theta(\{J > N^\alpha\}) \le Ce^{-C_0 N^{\alpha-\beta}},$$

which means that $\mathbb{E}^2_\theta(J) \le N^\alpha + N\mathbb{P}^2_\theta(\{J > N^\alpha\}) \le CN^\alpha$. In view of (3.11), this proves (3.2), provided we have chosen $\delta$ sufficiently small that $\varrho + \alpha < 1 - \delta$. □

Our program is thus completed once we prove (3.7). To this end, an intermediate result is needed.

LEMMA 3.4. *Given $R > 0$, take two points $a_R$ and $b_R$ such that $\|a_R - b_R\| = R$. For each $\epsilon \in (0,1]$, consider two walks starting at $a_R$ and $b_R$, respectively, and define the stopping time $\tau_{\epsilon,R}$ as the first time $n$ such that*

$$\|X_n - Y_n\| \le \frac{R\epsilon}{2} \quad \text{or} \quad \|X_n - Y_n\| \ge 2R.$$

*There then exist $R_\epsilon \in \mathbb{R}_+$ and $C_2 > 0$ such that for each $R \ge R_\epsilon$ and each $\theta \in \Theta$,*

$$\begin{aligned}&\mathbb{P}^2_\theta(\{\|X_{\tau_{\epsilon,R}} - Y_{\tau_{\epsilon,R}}\| \ge 2R\}) \ge \tfrac{1}{2} - C_2\epsilon, \\ &\mathbb{P}^2_\theta(\{\|X_{\tau_{1,R}} - Y_{\tau_{1,R}}\| \ge 2R\}) \ge \tfrac{1}{4}.\end{aligned}$$

PROOF. Of course, the estimate in the statement is essentially sharp only in one dimension. If $d > 1$, then the probability is actually close to one. Yet the above estimate suffices for our purposes. So, in the higher-dimensional case, we will control only one coordinate, whereby we obtain the same estimate as in one dimension.

We decompose $(\widehat{X}_n, \widehat{Y}_n)$ in the same way as we decomposed $\widehat{X}_n$ in (2.14). Observe that $\mathbb{E}(\widehat{\Delta}^X_{t+1} | \mathcal{F}^{0,X}_t) = \mathbb{E}(\widehat{\Delta}^X_{t+1} | \mathcal{F}^{0,XY}_t)$. Define

$$M^{XY}_n := (\widehat{X}_n, \widehat{Y}_n) - (h(\omega^X_0), h(\omega^Y_0)) + (h(\omega^X_n), h(\omega^Y_n)).$$

As before, $M^{XY}_n$ is a bounded martingale with respect to the filtration $\mathcal{F}^{0,XY}_n$, while the remainder is a bounded boundary term. Since $\|a_R - b_R\| = R$, it follows that there exists a unit vector $v_R$ such that $\langle v_R, X_0 - Y_0\rangle = R$. If we now define a new stopping time $\tau^*_R$ as the first time $n$ for which

$$\langle v_R, X_n - Y_n\rangle \le \frac{\epsilon R}{2} \quad \text{or} \quad \langle v_R, X_n - Y_n\rangle \ge 2R,$$

it then follows that

$$p := \mathbb{P}^2_\theta(\{\langle v_R, X_{\tau^*_R} - Y_{\tau^*_R}\rangle \geq 2R\}) \leq \mathbb{P}^2_\theta(\{\|X_{\tau_R} - Y_{\tau_R}\| \geq 2R\}).$$

Note that $\mathbb{E}^2_\theta(\tau^*_R) < \infty$. Indeed, by the ellipticity Assumption (A5), $\inf_{\theta,n}\mathbb{P}^2_\theta(\{\sup_{0\leq i\leq 2R}|\langle v_R, X_{n+i} - Y_{n+i}\rangle| \geq 2R\}) \geq \beta > 0$. Hence, it follows that $\mathbb{P}^2_\theta(\{\tau^*_R \geq n\}) \leq C(1-\beta)^{n/(2R)}$. Thus,

$$\begin{aligned}\mathbb{E}^2_\theta(\langle v_R, X_{\tau^*_R} - Y_{\tau^*_R}\rangle) &= \mathbb{E}^2_\theta(\langle (v_R, -v_R), M^{XY}_{\tau^*_R}\rangle) + \mathcal{O}(1)\\ &= \mathbb{E}^2_\theta(\langle (v_R, -v_R), M^{XY}_0\rangle) + \mathcal{O}(1) = R + \mathcal{O}(1),\end{aligned}$$

while, on the other hand,

$$\mathbb{E}^2_\theta(\langle v_R, X_{\tau^*_R} - Y_{\tau^*_R}\rangle) \leq 2Rp + \frac{\epsilon R}{2}(1-p).$$

The above two equations readily imply $p \geq \frac{1}{2} - \frac{\epsilon}{4} - CR^{-1}$, which is what we wanted. The second inequality follows similarly. □

We can conclude by proving (3.7).

Proof of Lemma 3.3. Let $X_0 = a$ and $Y_0 = b$, with $\|a - b\| \geq L_N$ and $\kappa \in (\frac{1}{2}, 1)$. Using the ellipticity Assumption (A5) for the first $L_N$ steps, the second estimate of Lemma 3.4 for the next $\ln_2 \epsilon^{-1}$ steps and, finally, the first estimate of that lemma for another $\log_2 N^\kappa$ steps, we obtain

$$\begin{aligned}&\mathbb{P}^2_\theta(\|X_j - Y_j\| \text{ reaches } N^\kappa L_N \text{ before } L_N)\\ &\qquad \geq \gamma^{L_N}\epsilon^{-2}(\tfrac{1}{2} - C_2\epsilon)^{\log_2 N^\kappa} \geq \epsilon^{-2}N^{-\kappa - C_3\epsilon - A\ln\gamma^{-1}}.\end{aligned}$$

In other words, there is a polynomially small probability of making an excursion of size $N^\kappa L_N$ before returning to a distance $L_N$. On the other hand, once we have such a big excursion, Hoeffding's inequality (see, e.g., [12]) implies that it will take more than $N$ steps to come back, indeed

$$\begin{aligned}&\mathbb{P}^2_\theta\left(\left\{\inf_{1\leq j\leq N}\|X_j - Y_j\| \leq L_N\right\} \,\Big|\, \|X_0 - Y_0\| \geq N^\kappa L_N\right)\\ &\qquad \leq 2\sup_\theta \mathbb{P}_\theta\left(\left\{\sup_{1\leq j\leq N}\|\widehat{X}_j\| > N^\kappa L_N/3\right\} \,\Big|\, \widehat{X}_0 = 0\right) \leq CNe^{-CN^{2\kappa-1}}.\end{aligned}$$

The last two inequalities imply (3.7) with $\beta = \kappa + C_3\epsilon + A\ln\gamma^{-1}$ and $C_0 = C\epsilon^{-2}$, provided we choose $A > 0$ and $\epsilon > 0$ sufficiently small that

$$\kappa + C_3\epsilon + A\ln\gamma^{-1} < 1. \tag{3.12}$$

This proves Lemma 3.3. □

**4. Proofs: the environment.** In this section, we prove Theorem 2.

REMARK 4.1. The reader is alerted to the fact that the following arguments are more of a functional analytic than of a probabilistic nature. In particular, we will work with the Banach space $\mathcal{M}(\Theta)$ of complex-valued measures, rather than just with probability measures.

4.1. *The environment*: *time mixing.* Following [15], we will lift the dynamics to a rather abstract space and prove that such a lift enjoys a spectral gap. This will imply the desired results.

In fact, we want to lift the dynamics to the space $\mathcal{B} := \overline{\mathbb{C} \otimes \bigotimes_{p\in\mathbb{Z}^d} \mathcal{M}_p(\Theta)}^{\|\cdot\|}$, where $\mathcal{M}_p(\Theta) := \{\mu \in \mathcal{M}(\Theta) : \mu(\varphi) = 0\ \forall \varphi \in \mathcal{C}^0(\Theta)$ that do not depend on $\theta_p\}$ and the closure is taken with respect to the norm

$$\|\bar{\boldsymbol{\mu}}\| := \sup\{|c_\mu|, |\mu_p| : p \in \mathbb{Z}^d\}.$$

Here, we use the notational convention that an element $\bar{\boldsymbol{\mu}} \in \mathcal{B}$ has components $c_\mu \in \mathbb{C}$ and $\bar{\mu} := (\mu_p)_p$ with $\mu_p \in \mathcal{M}_p(\Theta)$ and, for each complex-valued measure $\mu \in \mathcal{M}(\Theta)$, $|\mu|$ is the total variation of $\mu$. For example, if $\nu_p, \nu'_p$ ($p \in \mathbb{Z}^d$) are probability measures on $\Theta$ such that $\nu_p = \nu'_p$ for all $p \neq q$ and we set $\nu := \bigotimes_{p\in\mathbb{Z}^d} \nu_p$, $\nu' := \bigotimes_{p\in\mathbb{Z}^d} \nu'_p$, then $\nu - \nu' \in \mathcal{M}_q(\Theta)$.

To define such a lift, we first need to define a map $\Psi : \mathcal{M}(\Theta) \to \mathcal{B}$ and a projection $\mathrm{Pr} : \mathcal{B} \to \mathcal{M}(\Theta)$ that allow us to transfer objects between the two settings.

The choice of the first map is quite arbitrary; we will fix a convenient one. Consider a strict total ordering $\prec$ of $\mathbb{Z}^d$ such that $0 \prec p$ for each $p \in \mathbb{Z}^d \setminus \{0\}$, and the set $\{q : q \prec p\}$ contains the box (centered at zero) of size $C_4\|p\|$ and is contained in the box of size $C_5\|p\|$. For example, one can start from zero and spiral out on larger and larger cubical shells. Let $q_+$ be the successor of $q$ (i.e., $q \prec q_+$ and there are no $q' \in \mathbb{Z}^d$ such that $q \prec q' \prec q_+$).

Let $m$ be an arbitrary probability measure on $\Theta$, fixed once and for all. For each $q \in \mathbb{Z}^d$, we can then consider the $\sigma$-algebra $\mathcal{F}_q$ determined by all the variables $\omega^{q'}$ with $q' \prec q$, hence $\mathcal{F}_0$ is the trivial $\sigma$-algebra. Call $x^{\prec p}$ and $x^{\succ p}$ the set of coordinates with indices smaller (resp. larger) than $p$.

Next, for each $q \in \mathbb{Z}^d$, define the operator $J_q : \mathcal{C}^0(\Theta) \to \mathcal{C}^0(\Theta)$ and its dual $J'_q : \mathcal{M}(\Theta) \to \mathcal{M}(\Theta)$ by

$$J_q f := m^q(f) - m^{q+}(f), \qquad J'_q\mu(f) := \mu(J_q f), \tag{4.1}$$

where $m^q$ is the marginal of $m$ on the $\mathcal{F}_q$ variables. In other words, given a local function $f \in \mathcal{C}^0(\Theta)$ and a point $q \in \mathbb{Z}^d$, $m^q(f)$ depends only on $x^{\succeq q}$ variables, namely

$$m^q(f)(x^{\succeq q}) = \int_\Theta f(y^{\prec q}, x^{\succeq q}) m(dy^{\prec q}, dy^{\succeq q}).$$

For each local function $f$, we have

$$f = m(f) + \sum_{q\in\mathbb{Z}^d} J_q(f).$$

Note that as $f$ is local, there exists a box $\Lambda \subset \Theta$ such that $f$ depends only on the variables $\{\omega_q : q \in \Lambda\}$, but this means that the sum consists of only finitely many terms. Accordingly, for each $\mu \in \mathcal{M}(\Theta)$, we define $\hat{\mu} := \mu - \mu(1)m$ and for each $f \in \mathcal{C}^0$ and $q \in \mathbb{Z}^d$, we consider $J'_q\hat{\mu}(f) = \hat{\mu}(J_q f)$. Then, $J'_q\hat{\mu} \in \mathcal{M}_q(\Theta)$ and we can define the lift

$$\Psi(\mu) := (\mu(1), (J'_q\hat{\mu})_q) \in \mathcal{B}.$$

Remark 4.2. If one chooses $m := \bigotimes_{q\in\mathbb{Z}^d} m_*$, that is, a product measure, where $m_*$ is an arbitrary probability measure on $I$, then a direct computation using definition (4.1) yields $J_q f = \mathbb{E}_m(f \mid \mathcal{F}_q^c) - \mathbb{E}_m(f \mid \mathcal{F}_{q+}^c)$, where $\mathcal{F}_q^c$ is the $\sigma$-algebra determined by the $\omega_p$ with $p \succeq q$ and $\mathbb{E}_m$ is the expectation with respect to $m$. Moreover, $J'_q m = 0$ for all $q \in \mathbb{Z}^d$, hence $J'_q\hat{\mu} = J'_q\mu$. The reason to allow nonproduct measures in spite of the slightly more complex definitions is their usefulness in Section 4.5.

On the other hand, for each $\bar{\boldsymbol{\mu}} = (c_\mu, (\mu_p)_p) \in \mathcal{B}$ and local function $f$, we can define

$$\Pr \bar{\boldsymbol{\mu}}(f) := c_\mu m(f) + \sum_{p\in\mathbb{Z}^d} \mu_p(f). \tag{4.2}$$

Remark 4.3. Although $\Pr \bar{\boldsymbol{\mu}}(f)$ is well defined on each local function, $\Pr \bar{\boldsymbol{\mu}}$ is not necessarily a measure. Yet, $\Pr \Psi\mu = \mu$ for all $\mu \in \mathcal{M}(\Theta)$ since for each local function $f$, $m(\sum_q J_q f) = m(f - m(f)) = 0$. There thus exists a subset $\mathcal{B}_m \subset \mathcal{B}$ containing $\Psi(\mathcal{M}(\Theta))$ such that each element of $\Pr \mathcal{B}_m$ gives rise to a bounded linear functional on the space of local functions, hence uniquely identifies a measure. In other words, for each $\bar{\boldsymbol{\mu}} = (c_\mu, (\mu_p)_p) \in \mathcal{B}_m$, $c_\mu m + \sum_{p\in\mathbb{Z}^d} \mu_p$ converges weakly to a measure that we call $\Pr \bar{\boldsymbol{\mu}}$. (Note that here the order of the series may matter; we tacitly assume that the order is the one given by the relation $\prec$.)

Now that we know how to lift measures, we must describe how to lift the dynamics. For each local function $f$ and $p \in \mathbb{Z}^d$, we decompose

$$\mathcal{K} f = \mathcal{K}_p f + \sum_{q\in\mathbb{Z}^d\setminus\{p\}} \mathcal{K}_{p,q} f \tag{4.3}$$

with operators $\mathcal{K}_p$ and $\mathcal{K}_{p,q}$ defined as follows: for each $q \in \mathbb{Z}^d$, setting $q' \prec_p q :\Longleftrightarrow q' - p \prec q - p$,

$$
\begin{aligned}
(\mathcal{K}_{p,q} f)(\omega) := \int_\Theta K(\tau^p \omega, dy^p) &\prod_{0 \prec_p q' \prec_p q} K(\tau^{q'} \omega_{(p)}, dy^{q'}) \\
&\times [K(\tau^q \omega, dy^q) - K(\tau^q \omega_{(p)}, dy^q)] \\
&\times \prod_{q' \succ_p q} K(\tau^{q'} \omega, dy^{q'}) f(y),
\end{aligned} \tag{4.4}
$$

$$
(\mathcal{K}_p f)(\omega) := \int_\Theta K(\tau^p \omega, dy^p) \prod_{q' \in \mathbb{Z}^d \setminus \{p\}} K(\tau^{q'} \omega_{(p)}, dy^{q'}) f(y),
$$

where $\omega_{(p)}^{q'} = \omega^{q'}$ for each $q' \neq p$, while $\omega_{(p)}^p = a$ for some fixed $a \in I$. It is easy to see that the series in (4.3) converges due to Assumptions (A7) and (A9)(a). The fundamental fact of the above decomposition is that if $f$ does not depend on $\omega^q$, then $\mathcal{K}_{p,q} f = 0$. Accordingly, for each $\mu \in \mathcal{M}(\Theta)$, $\mathcal{K}'_{p,q}\mu \in \mathcal{M}_q(\Theta)$. In addition, if $\mu \in \mathcal{M}_p(\Theta)$, then $\mathcal{K}'_p \mu \in \mathcal{M}_p(\Theta)$ since if $f(y)$ does not depend on $y^p$, then $\mathcal{K}_p f$ also does not. Define $\bar{\boldsymbol{\alpha}} = (1, \bar{\alpha}) = (1, (\alpha_p)_p)$ by

$$
\bar{\boldsymbol{\alpha}} := \Psi(\mathcal{K}' m). \tag{4.5}
$$

Also, let

$$
\mathcal{K}' \mu_p = \sum_{q \in \mathbb{Z}^d \setminus \{p\}} \mathcal{K}'_{p,q} \mu_p + \mathcal{K}'_p \mu_p. \tag{4.6}
$$

Finally, (4.6) suggests that we define the operator $\overline{\mathcal{K}} : \mathcal{B} \to \mathcal{B}$ by

$$
\overline{\mathcal{K}} \bar{\boldsymbol{\mu}} := (c_\mu, c_\mu \bar{\alpha} + A\bar{\mu}) := \left( c_\mu, \left( c_\mu \alpha_q + \mathcal{K}'_q \mu_q + \sum_{p \in \mathbb{Z}^d \setminus \{q\}} \mathcal{K}'_{p,q} \mu_p \right)_q \right). \tag{4.7}
$$

For all $\mu \in \mathcal{M}(\Theta)$, local functions $f$ and $n \in \mathbb{N}$, we have

$$
\Pr \overline{\mathcal{K}}^n \Psi \mu(f) = \mu(\mathcal{K}^n f).
$$

Thus, $\overline{\mathcal{K}}(\mathcal{B}_m) \subseteq \mathcal{B}_m$ and the dynamics of $\overline{\mathcal{K}}$ covers the original one.

LEMMA 4.4. *The hypotheses (A6), (A7) and* (A9)(a) *on the Markov process imply* $\|A\| \leq \sum_q d_q = \eta_0 < 1$.

PROOF. Let $\nu \in \mathcal{M}_q(\Theta)$ and $f \in \mathcal{C}^0$. Define $\tilde{f}$ by $\mathcal{K}_q f(\omega) = \int K(\tau^q \omega, dy^q) \tilde{f}(y^q, \omega)$. Then, $|\tilde{f}|_\infty \leq |f|_\infty$ and $\tilde{f}$ does not depend on $\omega^q$. Now, by Assumption (A6), it follows that, varying $\omega^q$, $\int K(\tau^q \omega, dy^q) \tilde{f}(y^q, \omega)$ changes by at most $2d_0$. Hence, for each $(\omega^p)_{p \neq q}$, there must exist $a', a'' \in I$ and $t \in [0,1]$ such that defining $\omega'$ and $\omega''$ as those configurations which are obtained from a configuration $\omega$ replacing $\omega^0$ by $a'$ or $a''$, respectively, and defining $\tilde{K}(\omega,\cdot) = tK(\omega',\cdot) + (1-t)K(\omega'',\cdot)$, the following holds:

$$\left|\int K(\tau^q\omega, dy^q)\tilde{f}(y^q,\omega) - \int \tilde{K}(\tau^q\omega, dy^q)\tilde{f}(y^q,\omega)\right| \leq d_0|f|_\infty.$$

Thus,

$$\begin{aligned}|\nu(\mathcal{K}_q f)| &= \left|\nu\left(\int K(\tau^q\omega, dy^q)\tilde{f}(y^q,\omega) - \int \tilde{K}(\tau^q\omega_{(q)}, dy^q)\tilde{f}(y^q,\omega)\right)\right| \\ &\leq d_0|\nu||f|_\infty.\end{aligned}$$

While, by Assumption (A7), for all $p \neq q$ and all $\bar{\mu}$, we have $|\mathcal{K}'_{p,q}\mu_p| \leq d_{p-q}|\mu_p|$. □

Hence, the fixed point equation $\overline{\mathcal{K}}(1,\bar{\mu}) = (1,\bar{\mu})$ has the unique solution $(1,(\mathbb{1}-A)^{-1}\bar{\alpha})$, which can easily be seen to project down to a stationary probability measure $\mu_e$ on $\Theta$. Indeed, given any probability measure $\nu$, the set $\{\mathcal{K}^{n\prime}\nu\}$ will have weak accumulation points. On the other hand, calling $(1,\bar{\nu})$ the lift of $\nu$, we have that $\overline{\mathcal{K}}^n(1,\bar{\nu}) = (1, A^n\bar{\nu} + \sum_{k=0}^{n-1} A^k\bar{\alpha})$ is a lift of $\mathcal{K}^{n\prime}\nu$. Hence, the measures $\mathcal{K}^{n\prime}\nu$ must agree, on local functions, with the projections of the $\overline{\mathcal{K}}^n(1,\bar{\nu})$ and it follows that

$$\bar{\boldsymbol{\mu}} := \lim_{n\to\infty} \overline{\mathcal{K}}^n(1,\bar{\nu}) = (1,(\mathbb{1}-A)^{-1}\bar{\alpha}) \tag{4.8}$$

projects to a unique invariant probability measure $\mu_e$ which is the weak limit of the sequence $(\mathcal{K}^{n\prime}\nu)$.

Finally, the operator $\overline{\mathcal{K}}$ has a spectral gap, which implies exponential time-mixing of this invariant measure, that is, the analogue of property (A2) for the Markov evolution of the environment and the measure $\mu_e$. In particular, we have proven Assumption (A0) relative to $\mu_e$. [The fact that $\mu_e$ cannot be supported on the translation invariant configurations is a direct consequence of Assumption (A6); see also Lemma 4.4.]

REMARK 4.5. Note that the above argument would hold verbatim for more general, site-dependent, kernels $K_q(\tau^q\theta, dy^q)$ [instead of $K(\tau^q\theta, dy^q)$], the only difference being the loss of the translation invariance of $\mu_e$.

4.2. *The environment*: *space mixing.* For property (A3), we need an extra argument. Given a function $\varphi \in \mathcal{C}^0(\Theta)$, let us call $\Lambda(\varphi) \subset \mathbb{Z}^d$ the set of variables on which $\varphi$ depends. That is, $\Lambda(\varphi)$ is the smallest subset of $\mathbb{Z}^d$ such that $\varphi(x) = \varphi(y)$ whenever $x, y \in \Theta$ with $x^p = y^p$ for all $p \in \Lambda(\varphi)$. Also, let us call $\mathcal{C}^0_{\text{loc}}$ the set of continuous local functions [i.e., is functions for which $\Lambda(\varphi)$ is a finite set]. Finally, for any two functions $\varphi, \psi \in \mathcal{C}^0(\Theta)$, let $\rho_{\varphi,\psi} := \inf\{\|x-y\| : x \in \Lambda(\varphi), y \in \Lambda(\psi)\}$, the distance between the sets of dependence.

For each $\varphi \in \mathcal{C}^0_{\text{loc}}$ and $\psi \in \mathcal{C}^0$ such that $\rho_{\varphi\psi} = L$ and $\Lambda(\varphi)$ is contained in a box of size $l$, we want to estimate $\mu_e(\varphi\psi) - \mu_e(\varphi)\mu_e(\psi)$. Let us define $\Lambda_L(\varphi) := \{q \in \mathbb{Z}^d : \inf_{p\in\Lambda(\varphi)} \|p-q\| < L\}$. Clearly, $\Lambda_L(\varphi) \cap \Lambda(\psi) = \varnothing$. Moreover, $\Lambda_L(\varphi)$ contains at most $(l+L)^d$ sites. Next, for each $n \in \mathbb{N}$, let $r := L/2n$. Our main idea is to modify the kernels so as to define a new process with coupling range less than $r$ in $\Lambda_L(\varphi)$. To do so, we define, for each $q \in \Lambda_L(\varphi)$, the cutoff kernels

$$K_{r,q}(\omega, dy) := \begin{cases} K(\omega_{(r,q)}, dy), & \text{if } q \in \Lambda_L(\varphi), \\ K(\tilde{\omega}_{(r,q)}, dy), & \text{if } q \notin \Lambda_L(\varphi), \end{cases}$$

where, for some fixed $b \in I$,

$$\omega^p_{(r,q)} = \begin{cases} \omega^p, & \text{if } \|p-q\| < r, \\ b, & \text{if } \|p-q\| \ge r, \end{cases}$$

$$\tilde{\omega}^p_{(r,q)} = \begin{cases} \omega^p, & \text{if } \|p-q\| < r \text{ or } p \notin \Lambda_L(\varphi), \\ b, & \text{if } \|p-q\| \ge r \text{ and } p \in \Lambda_L(\varphi). \end{cases}$$

We can then use the above kernels to define the operator ${}^r\mathcal{K}$ as in formula (2.23) (see also Remark 4.5). Note that such an operator is close to the original one; indeed, for each $\phi \in \mathcal{C}^0(\Theta)$,

$$\begin{aligned} \|\mathcal{K}\phi - {}^r\mathcal{K}\phi\|_\infty &\le \sum_{q\in\Lambda_L(\varphi)} \sum_{\|z\|\ge r} d_z \|\phi\|_\infty + \sum_{q\notin\Lambda_L(\varphi)} \sum_{\substack{p\in\Lambda_L(\varphi) \\ \|p-q\|\ge r}} d_{p-q}\|\phi\|_\infty \\ &\le C(l+L)^d r^{-\xi'} \|\phi\|_\infty, \end{aligned} \tag{4.9}$$

by the long-range Assumption (A8). We can thus write

$$\begin{aligned} \mu_e(\varphi\psi) &= \mu_e(\mathcal{K}^n(\varphi\psi)) \\ &= \mu_e({}^r\mathcal{K}^n(\varphi\psi)) + \mathcal{O}(n(l+L)^d r^{-\xi'} \|\varphi\psi\|_\infty) \\ &= \mu_e(({}^r\mathcal{K}^n\varphi)({}^r\mathcal{K}^n\psi)) + \mathcal{O}(n(l+L)^d r^{-\xi'} \|\varphi\psi\|_\infty). \end{aligned} \tag{4.10}$$

At this point, it is natural to define measures $\nu_{n,r}(\phi) := \mu_e(\phi\, {}^r\mathcal{K}^n\psi)$. Note that the lift of such a measure to the space $\mathcal{B}$ is given by $\Psi(\nu_{n,r}) = (\mu_e({}^r\mathcal{K}^n\psi), \bar{\nu}_{n,r})$ with $\|\Psi(\nu_{n,r})\| \le 4\|\psi\|_\infty$. By the results of the previous section (and Remark 4.5) applied to the operator ${}^r\mathcal{K}$, it follows that there exists a measure $\mu_r$ such that ${}^r\mathcal{K}'\mu_r = \mu_r$ and, in addition,

$$|\nu_{n,r}({}^r\mathcal{K}^n\varphi) - \mu_r(\varphi)\mu_e({}^r\mathcal{K}^n\psi)| \le C\sigma^n l^d \|\varphi\|_\infty \|\psi\|_\infty.$$

Observe that due to the constructive nature of the proof in the preceding section, the constants $C$ and $\sigma$ do not depend on $r$.

The above estimates applied to the case $\psi = 1$ (i.e., to $\nu_{n,r} = \mu_e$) imply

$$|\mu_e(\varphi) - \mu_r(\varphi)| \le C\sigma^n l^d \|\varphi\|_\infty + Cn(l+L)^d r^{-\xi'}\|\varphi\|_\infty.$$

Choosing $n$ proportional to $\ln L$, the two last facts and (4.10) together yield, for $l \le L$,

$$\begin{aligned} &|\mu_e(\varphi\psi) - \mu_e(\varphi)\mu_e(\psi)| \\ &\qquad \le |\nu_{n,r}({}^r\mathcal{K}^n\varphi) - \mu_e(\varphi)\mu_e(\psi)| + C\frac{n(l+L)^d}{r^{\xi'}}\|\varphi\psi\|_\infty \\ &\qquad \le |\mu_e(\varphi)||\mu_e({}^r\mathcal{K}^n\psi) - \mu_e(\psi)| \\ &\qquad\quad + C(nL^d r^{-\xi'} + \sigma^n l^d)\|\varphi\|_\infty\|\psi\|_\infty \\ &\qquad \le C(L^{-\xi'+d}(\ln L)^{\xi'+1}\|\varphi\|_\infty\|\psi\|_\infty). \end{aligned} \tag{4.11}$$

Clearly, (4.11) implies the space-mixing property (A3), provided $\xi' > \xi + d$.

4.3. *The environment as seen from the particle.* The above construction can also be used to achieve our other goal—the study of the dynamics as seen from the particle. That is, we wish to study the operator [see (2.5)]

$$Sf(\omega) = \sum_{z\in\Lambda} \pi_z(\omega)\mathcal{K}(\tau^z f)(\omega). \tag{4.12}$$

We will use the same space as in the previous section, the only new difficulty being to define the covering dynamics. To this end, let us define $\mathbf{B}_z : \mathcal{B} \to \mathcal{B}$ by $\mathbf{B}_z\bar{\boldsymbol{\mu}} := (c_\mu, ((\mathbf{B}_z\bar{\boldsymbol{\mu}})_q)_q)$, where $(\mathbf{B}_z\bar{\boldsymbol{\mu}})_q(f) := \mu_{q+z}(\pi_z\tau^z f)$. We have $(\mathbf{B}_z\bar{\boldsymbol{\mu}})_q \in \mathcal{M}_q(\Theta)$, provided $q \notin \Lambda$, but we have no control of this kind for $q \in \Lambda$. Thus, although the operator $\sum_{z\in\Lambda}\overline{\mathcal{K}}\mathbf{B}_z$ would cover $S'$, it does not respect the structure of the space for the components in $\Lambda$. On the other hand, $(\mathbf{A}_z\bar{\boldsymbol{\mu}})_q(f) := \mu_{q+z}(a_z\tau^z f) \in \mathcal{M}_q(\Theta)$ for each $q \in \mathbb{Z}^d$. Accordingly, for each $q \in \Lambda$, we must deal only with the remainders $(\widehat{\mathbf{R}}_z\bar{\boldsymbol{\mu}})_q(f) := \mu_{q+z}(\hat{\pi}_z\tau^z f)$. We can again use our decomposition operators to write $(\widehat{\mathbf{R}}_z\bar{\boldsymbol{\mu}})_q(f) = \sum_{q'\in\mathbb{Z}^d}(\widehat{\mathbf{R}}_z\bar{\boldsymbol{\mu}})_q(J_{q'}f)$ and then redistribute the various terms to the appropriate components of the vector; indeed, $J'_{q'}(\widehat{\mathbf{R}}_z\bar{\boldsymbol{\mu}})_q \in \mathcal{M}_{q'}(\Theta)$. Finally, defining, in analogy with the previous case, $\bar{\zeta}$ by $\Psi(S'\mu_e) = (1, \bar{\zeta})$, we can define the covering dynamics $\overline{\mathbf{S}}(c_\nu, \bar{\nu}) := (c_\nu, c_\nu\bar{\zeta} + \overline{S}\bar{\nu})$, where

$$(\overline{S}\bar{\nu})_q := \begin{cases} \displaystyle\sum_{z\in\Lambda}\left[(\overline{\mathcal{K}}\mathbf{B}_z(0,\bar{\nu}))_q + \sum_{q'\in\Lambda} J'_q(\overline{\mathcal{K}}\widehat{\mathbf{R}}_z(0,\bar{\nu}))_{q'}\right], & \text{if } q \notin \Lambda, \\ \displaystyle\sum_{z\in\Lambda}\left[(\overline{\mathcal{K}}\mathbf{A}_z(0,\bar{\nu}))_q + \sum_{q'\in\Lambda} J'_q(\overline{\mathcal{K}}\widehat{\mathbf{R}}_z(0,\bar{\nu}))_{q'}\right], & \text{if } q \in \Lambda. \end{cases} \tag{4.13}$$

Once more, one can easily check that for each $\nu \in \mathcal{M}(\Theta)$, the following holds:

$$\Pr[\overline{\mathbf{S}}^n \Psi\nu](f) = \nu(S^n f) \qquad \text{for all } n \in \mathbb{N},$$

that is, we have a proper covering of the original dynamics. Moreover, by (A9)(b),

$$\begin{aligned} \|\overline{S}\| &\le \|\overline{\mathcal{K}}\| \left(1 + (1 + 2|\Lambda|) \sum_{z\in\Lambda} |\hat{\pi}_z|_\infty\right) \\ &\le \eta_0(1 + (1 + 2|\Lambda|)D) = \eta_1 < 1. \end{aligned} \tag{4.14}$$

This proves properties (A0) for $\mu$ and (A2) by exactly the same arguments used in Section 4.1.

4.4. *Locality.* The idea for the verification of condition (A4) is, again, to approximate the true dynamics with one in which the interactions are cut off at a proper scale. More precisely, let $A', B'$ be two $L/2$-neighborhoods of $A, B$, respectively. Similarly to what we did in the previous section, we can kill all interactions in $A' \cup B'$ at a distance larger than $L/2s$. If we call $\mathbb{P}^0_\theta$ the distribution of the resulting Markov process started from the configuration $\theta$, then Assumption (A8) implies

$$\begin{aligned} &|\mathbb{P}^e_\theta(f(\theta_1, \dots, \theta_s)g(\theta_1, \dots, \theta_s)) - \mathbb{P}^0_\theta(f(\theta_1, \dots, \theta_s)g(\theta_1, \dots, \theta_s))| \\ &\quad \le Cs(M + L)^d s^{\xi'} L^{-\xi'} \|fg\|_\infty \le C_M L^{-(\xi'-d)} s^{\xi'+1} \|fg\|_\infty \end{aligned}$$

with $\xi' > d + 2$. As $\mathbb{P}^0_\theta(f(\theta_1, \dots, \theta_s)g(\theta_1, \dots, \theta_s)) = \mathbb{P}^0_\theta(f(\theta_1, \dots, \theta_s))\mathbb{P}^0_\theta(g(\theta_1, \dots, \theta_s))$, by construction, condition (A4) follows, provided we choose $\xi = \xi' - d > 4$ and $\tilde{\xi} > \xi' + 1$.

4.5. *Absolute continuity.* Since we aim to prove that $\mu$ is absolutely continuous with respect to $\mu_e$, it is natural to work in the smaller space $\mathcal{M}^e$ of measures on $\Theta$ which are absolutely continuous with respect to $\mu_e$. In addition, we now have a natural reference measure, so we want to choose the measure $m$, in the construction of the lift $\Psi$ defined at the beginning of Section 4.1, to be $\mu_e$. This implies $\Psi(\mu_e) = (1, 0) \in \mathcal{B}$. Clearly, all the results of the previous section apply with this choice of $m$.

First, note that the above is consistent.

LEMMA 4.6. *Both $\mathcal{K}'$ and $S'$ are well defined as operators from $\mathcal{M}^e$ to itself.*

PROOF. First consider the case of a measure $\nu \in \mathcal{M}^e$ such that $f := \frac{d\nu}{d\mu_e} \in L^\infty(\Theta, \mu_e)$. Then,

$$\mathcal{K}'(\nu)(\varphi) = \int f \cdot \mathcal{K}\varphi \, d\mu_e \le \|f\|_\infty \cdot \int \mathcal{K}|\varphi| \, d\mu_e = \|f\|_\infty \cdot \mu_e(|\varphi|)$$

for each bounded measurable $\varphi:\Theta\to\mathbb{R}$. Now, for each $\nu\in\mathcal{M}^e$, again calling $f$ the density, by monotone convergence,

$$\mathcal{K}'(\nu)(\varphi)=\sup_n \mathcal{K}'((f\wedge n)\mu_e)(\varphi)\le \sup_n n\mu_e(\varphi)$$

for each such $\varphi$. In particular, this holds for $\varphi=1_A$ being the indicator function of any $\mu_e$-null set. Hence, $\mathcal{K}'(\nu)(A)=0$ if $\mu_e(A)=0$. Similar arguments hold for $S'$ by using the translation invariance of $\mu_e$. $\square$

On the other hand, in the present generality, it is not guaranteed that $J_q$ maps $\mathcal{M}^e$ into itself, so we cannot check directly that the covering dynamics $\overline{\mathbf{S}}:\mathcal{B}\to\mathcal{B}$ preserves absolute continuity of the components.

The above considerations show that it may not be very convenient to work with the decomposition (4.13) to treat the problem of absolute continuity. Furthermore, the second sum on the left-hand side of (4.13) is of a highly nonlocal nature, which makes it very hard to control the sum of the components of the resulting vector. To overcome these difficulties, it is useful to decompose the operator $S$ into a convex combination of two operators—one representing a random walk with fixed transition probabilities and the other (small) one keeping the dependence on the environment. To do so, it suffices to write

$$\pi_z=(1-\kappa)c_z+\kappa\tilde{\pi}_z, \tag{4.15}$$

where $c_z:=(1-\kappa)^{-1}\max\{0,a_z-\|\hat{\pi}_z\|_\infty\}$, $\tilde{\pi}_z:=\kappa^{-1}(\pi_z-\max\{0,a_z-\|\hat{\pi}_z\|_\infty\})$ and $\kappa:=1-\sum_{z\in\Lambda}\max\{0,a_z-\|\hat{\pi}_z\|_\infty\}\le D$ is small by Assumption (A9)(c). Let us set $S_0f:=\sum_z c_z\tau^z\mathcal{K}f$ and $S_1:=\sum_z\tilde{\pi}_z\tau^z\mathcal{K}f$. Clearly $S=(1-\kappa)S_0+\kappa S_1$ and $c_z,\tilde{\pi}_z\ge 0$, $\sum_z c_z=\sum_z\tilde{\pi}_z=1$. Hence, $|S_0f|_\infty\le|f|_\infty$, $|S_1f|_\infty\le|f|_\infty$. It is then convenient to consider a Bernoulli process with probability $(1-\kappa,\kappa)$. For each $\sigma\in\{0,1\}^{\mathbb{N}}$, we let $\sigma^n$ be the first $n$ symbols of $\sigma$ and set $S_{\sigma^n}:=S_{\sigma_n}\cdots S_{\sigma_1}$. Using $E$ for the expectation with respect to the above process, we see that

$$S^nf=E(S_{\sigma^n}f). \tag{4.16}$$

The advantage of the representation (4.16) lies in the fact that $S_0$ can be conveniently treated by our covering space techniques, while the occurrences of $S_1$ are weighted by a small probability.

To be more precise, recall that, by the analogue of (4.8), $\mu=\lim_{N\to\infty}S^{N\prime}\mu_e$. Then, after setting $\zeta:=(S'-\mathbb{1})\mu_e$, we can write

$$\begin{gathered}S^{N\prime}\mu_e=\mu_e+\sum_{n=0}^{N-1}S^{n\prime}\zeta,\\ |S^{n\prime}\zeta(\varphi)|\le|E(\zeta(S_{\sigma^n}\varphi)\mathbb{1}_{\Sigma_n})|+E(\mathbb{1}_{\Sigma_n^c})\|\varphi\|_\infty,\end{gathered} \tag{4.17}$$

where $\Sigma_n$ is the set of $\sigma$ such that $\sigma^n$ contains a string of zeros with length greater than $t_n := -(1-\vartheta)[\ln(1-\kappa)]^{-1}\ln n$, $\vartheta \in (0,1)$. Clearly, $E(\mathbb{1}_{\Sigma_n^c}) \le e^{-C_\vartheta n^{\vartheta/2}}$. On $\Sigma_n$, let $m_\sigma$ denote the beginning of the first string of $t_n$ zeros. We have

$$\begin{aligned}|\zeta(S^n\varphi)| &\le \left|\sum_{q\in\mathbb{Z}^d} E(\zeta(S_{\sigma_n,\dots,\sigma_{m_\sigma+t_n}} J_q S_0^{t_n} S_{\sigma^{m_\sigma}}\varphi)1_{\Sigma_n})\right| + e^{-C_\vartheta n^{\vartheta/2}}\|\varphi\|_\infty \\ &\le E(|\Pr(\overline{\mathbf{S}}_0^{t_n}\Psi(S'_{\sigma_n,\dots,\sigma_{m_\sigma+t_n}}\zeta))(S_{\sigma^{m_\sigma}}\varphi)1_{\Sigma_n}|) + e^{-C_\vartheta n^{\vartheta/2}}\|\varphi\|_\infty,\end{aligned} \tag{4.18}$$

where we have introduced the covering operator $\overline{\mathbf{S}}_0$ defined by

$$\overline{\mathbf{S}}_0(c_\nu,(\bar\nu)_q) = (c_\nu,(\bar S_0\bar\nu)_q) := \left(c_\nu, \sum_{z\in\Lambda}(\overline{\mathcal{K}}\mathbf{C}_z(0,\bar\nu))_q\right) \tag{4.19}$$

and, contrary to (4.13), we have defined $(\mathbf{C}_z(0,\bar\nu))_q(f) := \nu_{q+z}(c_z\tau^z f)$. A direct computation shows that $\overline{\mathbf{S}}_0^n$ covers $S_0^n$, hence formula (4.18).

Note that the summands in the first line of (4.17) are absolutely continuous measures with respect to $\mu_e$, by Lemma 4.6. Hence, the total variation of such measures is the $L^1(\Theta,\mu_e)$-norm of their density.

Unlike $\bar S$, the operator $\bar S_0$ is reasonably local. To make precise such a locality, we introduce the norm $\|\bar\nu\|_1 := \sum_{p\in\mathbb{Z}^d}|\nu_p|$ and define $\mathcal{B}_1 := \{(0,\bar\nu)\in\mathcal{B} : \|\bar\nu\|_1 < \infty\}$.

LEMMA 4.7. *For each* $(0,\bar\nu)\in\mathcal{B}_1$,

$$\|\overline{S}_0\bar\nu\|_1 \le \eta_0\|\bar\nu\|_1.$$

PROOF. Let $(0,\bar\nu)\in\mathcal{B}_1$. Then,

$$\begin{aligned}\sum_{p\in\mathbb{Z}^d}(\overline{S}_0\bar\nu)_p| &\le \sum_{z\in\Lambda}c_z\sum_{p\in\mathbb{Z}^d}\left[|\mathcal{K}'_p\tau^z\nu_{p+z}| + \sum_{q\in\mathbb{Z}^d\setminus\{p\}}|\mathcal{K}'_{q,p}\tau^z\nu_{q+z}|\right] \\ &\le \sum_{z\in\Lambda}c_z\sum_{p\in\mathbb{Z}^d}\left[d_0|\nu_{p+z}| + \sum_{q\in\mathbb{Z}^d\setminus\{p\}}d_{q-p}|\nu_{q+z}|\right] \\ &\le d_0\|\bar\nu\|_1 + \sum_{q\in\mathbb{Z}^d}|\nu_q|\sum_{z\in\Lambda}c_z\sum_{v\in\mathbb{Z}^d\setminus\{0\}}d_v \\ &= \left(d_0 + \sum_{v\ne 0}d_v\right)\|\bar\nu\|_1.\end{aligned}$$

□

Next, we verify that the above norm is relevant to the problem at hand.

LEMMA 4.8. *For each $n \in \mathbb{N}$ and $\sigma \in \{0,1\}^{\mathbb{N}}$, we have*

$$\|\Psi(S'_{\sigma^n}\zeta)\|_1 = \sum_{q\in\mathbb{Z}^d} |J'_q S'_{\sigma^n}\zeta| \le Cn^{((\xi'+1)/(\xi'-d))d}.$$

PROOF. We start by studying $J'_q S'_{\sigma^n}\mu_e$. We will also write $S$ for $S_\sigma$, since the computation is exactly the same. Given $q \in \mathbb{Z}^d$, we can change the kernel $K$ for all points $\|p\| \le \frac{C_4}{2}\|q\|$ to have only interactions of range $C_4\|q\|(8n)^{-1}$. (The constant $C_4$ is defined at the beginning of Section 4.1.)

We call $\mathcal{K}^{(q)}$ and $S^{(q)}$ the resulting Markov operators for the process of the environment and the environment as seen from the particle, respectively. Clearly, by Assumption (A8),

$$\|(S^{(q)} - S)\varphi\|_\infty \le C\|q\|^{-\xi'+d}n^{\xi'}\|\varphi\|_\infty.$$

Hence, for all $m \le n$,

$$\begin{aligned}\|((S^{(q)})^m - S^m)\varphi\|_\infty &\le \sum_{k=0}^{m-1}\|(S^{(q)})^k(S^{(q)} - S)S^{m-1-k}\varphi\|_\infty \\ &\le C\|q\|^{-\xi'+d}n^{\xi'}m\|\varphi\|_\infty\end{aligned} \tag{4.20}$$

and the same holds for $\mathcal{K}^{(q)}$ and $\mathcal{K}$. On the other hand, if $n \le \frac{C_4}{2}\|q\|$, then

$$\begin{aligned}(S^{(q)})^n J_q\varphi &= \sum_{z_1,\dots,z_n\in\Lambda} \pi_{z_1}\mathcal{K}^{(q)}\tau^{z_1}[\pi_{z_2}\mathcal{K}^{(q)}\cdots\mathcal{K}^{(q)}\tau^{z_{n-1}}[\pi_{z_n}\mathcal{K}^{(q)}\tau^{z_n}J_q\varphi]\cdots] \\ &= \sum_{z_1,\dots,z_n\in\Lambda} \pi_{z_1}\mathcal{K}^{(q)}\tau^{z_1}[\pi_{z_2}\mathcal{K}^{(q)}\cdots\mathcal{K}^{(q)}\tau^{z_{n-1}}[\pi_{z_n}]\cdots] \\ &\qquad\times\mathcal{K}^{(q)}\tau^{z_1}[\mathcal{K}^{(q)}\cdots\tau^{z_{n-1}}[\mathcal{K}^{(q)}\tau^{z_n}J_q\varphi]\cdots] \\ &= \sum_{z_1,\dots,z_n\in\Lambda} \pi_{z_1}\mathcal{K}^{(q)}\tau^{z_1}[\pi_{z_2}\mathcal{K}^{(q)}\cdots\mathcal{K}^{(q)}\tau^{z_{n-1}}[\pi_{z_n}]\cdots] \\ &\qquad\times\tau^{z_1+\cdots+z_n}\mathcal{K}^n J_q\varphi + \mathcal{O}(\|q\|^{-\xi'+d}n^{\xi'+1}\|\varphi\|_\infty),\end{aligned}$$

where we have used the fact that if two functions $f, g$ have disjoint support, then $\mathcal{K}^{(q)}(fg) = \mathcal{K}^{(q)}f\cdot\mathcal{K}^{(q)}g$, the same considerations as in (4.20) and the translation invariance of $\mathcal{K}$. Set $z^{n-k} := (z_{k+1},\dots,z_n)$ and

$$\psi_{z^n} := \tau^{-z_1-\cdots-z_n}[\pi_{z_1}\mathcal{K}^{(q)}\tau^{z_1}[\pi_{z_2}\mathcal{K}^{(q)}\cdots\mathcal{K}^{(q)}\tau^{z_{n-1}}[\pi_{z_n}]\cdots]].$$

Then,

$$\begin{aligned}\psi_{z^n} &= [\tau^{-z_1-\cdots-z_n}\pi_{z_1}]\mathcal{K}^{(q)}[\tau^{-z_2-\cdots-z_n}\pi_{z_2}]\mathcal{K}^{(q)}\cdots\mathcal{K}^{(q)}[\tau^{-z_n}\pi_{z_n}] \\ &= [\tau^{-z_1-\cdots-z_n}\pi_{z_1}]\mathcal{K}^{(q)}\psi_{z^{n-1}}\end{aligned}$$

since, when applied to a function $f$ supported in the box $\{p\in\mathbb{Z}^d : \|p\| \le \frac{C_4}{4}\|q\|\}$, the operator $\mathcal{K}^{(q)}$ is invariant under translation by $\tau^{-z_1-\cdots-z_n}$, provided $nC_1 \le \frac{C_4}{8}\|q\|$, that is $\mathcal{K}^{(q)}\tau^{-z_1-\cdots-z_n}f = \tau^{-z_1-\cdots-z_n}\mathcal{K}^{(q)}f$.

REMARK 4.9. Note that later on in the proof, we will apply $\mathcal{K}^{(q)}$ to functions that are supported in a box of size $\frac{C_4}{2}\|q\|$ centered at the origin. Hence, by construction, such functions never see the kernels in which the interaction is long range and, accordingly, we can modify $\mathcal{K}^{(q)}$ to have interactions of length $\|q\|$ in all of $\mathbb{Z}^d$ without any change in the above formulae. We will call $\widetilde{\mathcal{K}}^{(q)}$ the resulting translation invariant object.

We set $\Xi_n := \sum_{z_1,\dots,z_n}\psi_{z^n}$ and write the above estimate as

$$|\mu_e(S^n J_q\varphi) - \mu_e(\Xi_n\cdot(\mathcal{K}^{(q)})^n J_q\varphi)| \le C\|q\|^{-\xi'+d}n^{\xi'+1}\|\varphi\|_\infty. \tag{4.21}$$

As $J_q\varphi$ does not depend on points $\|p\|\le C_4\|q\|$, (4.11) yields

$$\begin{aligned}&|\mu_e(S^n J_q\varphi) - \mu_e(\Xi_n)\cdot\mu_e((\mathcal{K}^{(q)})^n J_q\varphi)|\\ &\qquad \le C\|q\|^{-\xi'+d}n^{\xi'+1}\|\varphi\|_\infty + C\|q\|^{-\xi'+d}(\ln\|q\|)^{\xi'+1}\|\varphi\|_\infty\|\Xi_n\|_\infty\end{aligned}$$

and as $\mu_e((\mathcal{K}^{(q)})^n J_q\varphi) = \mu_e(J_q\varphi) + \mathcal{O}(\|q\|^{-\xi'+d}n^{\xi'+1}\|\varphi\|_\infty)$, it follows that

$$\begin{aligned}&|\mu_e(S^n J_q\varphi) - \mu_e(J_q\varphi)|\\ &\qquad \le |\mu_e(\Xi_n - 1)|\|\varphi\|_\infty + C\|q\|^{-\xi'+d}n^{\xi'+1}\|\varphi\|_\infty(1+\|\Xi_n\|_\infty)\\ &\qquad\quad + C\|q\|^{-\xi'+d}(\ln\|q\|)^{\xi'+1}\|\varphi\|_\infty\|\Xi_n\|_\infty.\end{aligned}$$

At this point, we can estimate the real objects of interest:

$$\begin{aligned}|\zeta(S^n J_q\varphi)| \le \sum_{i=0}^{1}|\mu_e(\Xi_{n+i}-1)|\|\varphi\|_\infty + C\|q\|^{-\xi'+d}n^{\xi'+1}\|\varphi\|_\infty\|\Xi_n\|_\infty\\ + C\|q\|^{-\xi'+d}(\ln\|q\|)^{\xi'+1}\|\varphi\|_\infty\|\Xi_n\|_\infty.\end{aligned} \tag{4.22}$$

To conclude, we must estimate $|\mu_e(\Xi_n - 1)|$ and the norm $\|\Xi_n\|_\infty$ on the right-hand side of the above equation. To do so, it is convenient to consider the operators $\mathcal{K}_g\varphi := g\widetilde{\mathcal{K}}^{(q)}\varphi$, $\hat{\mathcal{K}}_g\varphi := g\widetilde{\mathcal{K}}^{(q)}\varphi - g\mu_e^{(q)}(\varphi)$, where $\mu_e^{(q)}$ is the unique invariant measure of the operator $\widetilde{\mathcal{K}}^{(q)}$. The proof of the existence of such measures is exactly the same as the proof of the existence of $\mu_e$ or $\mu_r$.

Note that $\hat{\mathcal{K}}_g 1 = 0$. If we set $\pi_{z^{n-k}} := \tau^{-z_{k+1}-\cdots-z_n}\pi_{z_{k+1}}$, then

$$\begin{aligned}\Xi_n &= \sum_{z^n}\mathcal{K}_{\pi_{z^n}}\mathcal{K}_{\pi_{z^{n-1}}}\cdots\mathcal{K}_{\pi_{z^1}}1\\ &= \sum_{z^n}\hat{\mathcal{K}}_{\pi_{z^n}}\mathcal{K}_{\pi_{z^{n-1}}}\cdots\mathcal{K}_{\pi_{z^1}}1 + \pi_{z^n}\mu_e^{(q)}(\mathcal{K}_{\pi_{z^{n-1}}}\cdots\mathcal{K}_{\pi_{z^1}}1).\end{aligned} \tag{4.23}$$

The translation invariance discussed above implies

$$\sum_{z^n} |\pi_{z^n} \mu_e^{(q)}(\mathcal{K}_{\pi_{z^{n-1}}} \cdots \mathcal{K}_{\pi_{z^1}} 1)| \le \left(1 + \sum_z \|\pi_z\|_\infty\right) = 1 + D. \tag{4.24}$$

To conclude, we need to understand the properties of compositions of the operators $\hat{\mathcal{K}}_{\pi_z^j}$. To do so, we again use our covering spaces, as we did for the operators $\mathcal{K}, S$. Indeed, given an operator $\hat{\mathcal{K}}_g$ with $g = a + \hat{g}$, $a \in \mathbb{R}$, and $\hat{g}$ supported in the box $p + \Lambda$ centered around $p$, we can define the covering dynamics $\overline{\mathcal{K}}\bar{G}$ of $\hat{\mathcal{K}}'_g$, where $\overline{\mathcal{K}}$ is defined as in (4.7), while $\bar{G}(c_\nu, \bar{\nu}) := (0, G\bar{\nu})$ is defined, in analogy with (4.13), by

$$(G\bar{\nu})_q := \begin{cases} g\nu_q + \sum_{p' \in p+\Lambda} J'_q(\hat{g}\nu_{p'}), & \text{if } q \notin p + \Lambda, \\ a\nu_q + \sum_{p' \in p+\Lambda} J'_q(\hat{g}\nu_{p'}), & \text{if } q \in p + \Lambda. \end{cases}$$

A direct computation shows that $\Pr(\overline{\mathcal{K}}\bar{G}_1 \cdots \overline{\mathcal{K}}\bar{G}_n \Psi\nu) = \hat{\mathcal{K}}'_{g_1} \cdots \hat{\mathcal{K}}'_{g_n}\nu$, that is, the above operators cover arbitrary products of the operators of interest. In addition, if $\phi \in \mathcal{C}^0_{\mathrm{loc}}$ depends only of $M$ variables, then

$$\|\hat{\mathcal{K}}_{g_n} \cdots \hat{\mathcal{K}}_{g_1} \phi\|_\infty \le \left[\prod_{k=1}^n (a + (1 + 2|\Lambda|)\|\hat{g}_k\|_\infty) \sum_q d_q\right] M \|\phi\|_\infty. \tag{4.25}$$

Iterating the procedure in (4.23), (4.24) and using (4.25) yields

$$\|\Xi_n\|_\infty \le C \sum_{k=0}^{n-1} \eta_1^k (1 + D) < \infty. \tag{4.26}$$

Next, observe that $\mu_e^{(q)}(\Xi_n) = 1$, which can easily be checked by a direct computation that uses the invariance of $\mu_e^{(q)}$ under translations and under the kernel $\tilde{\mathcal{K}}^{(q)}$. Hence, $|\mu_e(\Xi_n - 1)| \le |\mu_e(\Xi_n) - \mu_e^{(q)}(\Xi_n)|$ and as $\Xi_n$ depends on at most $C\|q\|^d$ variables, we can apply (4.8) and its analogue for the kernel $\tilde{\mathcal{K}}^{(q)}$ to conclude that

$$\begin{aligned} |\mu_e(\Xi_n) - \mu_e^{(q)}(\Xi_n)| &\le \|\mathcal{K}^n \Xi_n - (\tilde{\mathcal{K}}^{(q)})^n(\Xi_n)\|_\infty + C\|q\|^d \eta^n \|\Xi_n\|_\infty \\ &\le n \cdot C\|q\|^d \sum_{\|p\| \ge C_4\|q\|/8n} d_p + C\|q\|^d \eta^n \|\Xi_n\|_\infty \\ &\le C\|q\|^{d-\xi'} n^{1+\xi'} \cdot \|\Xi_n\|_\infty, \end{aligned} \tag{4.27}$$

where we used Assumption (A8) in the last step.

Combining (4.22), (4.26) and (4.27), we thus arrive at

$$|\zeta(S^n J_q \varphi)| \le C\|q\|^{d-\xi'} n^{1+\xi'} \cdot \|\varphi\|_\infty (1 + \|\Xi_n\|_\infty),$$

which, using the trivial bound $|J_q'(S')^n\nu| \le 2|\nu|$ for $\|q\| \le n^{(\xi'+1)/(\xi'-d)}$, implies the lemma, as $\xi' > 2d$ by Assumption (A8). □

Finally, applying Lemmas 4.8 and 4.7 to (4.18), we have

$$|\zeta(S^n\varphi)| \le C(\eta_0^{t_n} n^{((\xi'+1)/(\xi'-d))d} + e^{-C_\vartheta n^{\vartheta/2}})\|\varphi\|_\infty \le Cn^{-C_6}\|\varphi\|_\infty, \tag{4.28}$$

where $C_6 > 1$. Indeed, $\eta_0^{t_n} \le n^{-((1-\vartheta)\ln\eta_0^{-1})/(\ln(1-D)^{-1})}$, thus the claim holds true, provided $\frac{(1-\vartheta)\ln\eta_0^{-1}}{\ln(1-D)^{-1}} - \frac{\xi'+1}{\xi'-d}d > 1$ or, equivalently, $\eta_0 < (1-D)^{(\xi'(1+d))/((\xi'-d)(1-\vartheta))}$. By Assumption (A9)(c), we can always choose $\vartheta$ so that this inequality is satisfied.

Accordingly, the sum on the right-hand side of the first line of (4.17) is convergent, which implies that $\frac{d(S^n)'\mu_e}{d\mu_e}$ converges in $L^1(\Theta,\mu_e)$ to some function $h$, hence $\mu = h\mu_e$ is absolutely continuous with respect to $\mu_e$.

Finally, to prove equivalence, if $\mu_e$ is not absolutely continuous with respect to $\mu$, then there exists an invariant set $A$ such that $\mu(A) = 0$ but $\mu_e(A) > 0$. Accordingly,

$$0 = \mu(S^n\mathbb{1}_A) = \mu_e(hS^n\mathbb{1}_A) = \sum_{z_1,\dots,z_n\in\Lambda} \mu_e(h\pi_{z_1}\mathcal{K}\tau^{z_1}\pi_{z_2}\cdots\mathcal{K}\mathbb{1}_A),$$

which implies $\mu_e(h\mathcal{K}^n\tau^{\sum_{i=1}^n z_i}\mathbb{1}_A) = 0$ for each choice of $n$ and $z_i$ such that $\gamma_{z_i} \ne 0$. By ellipticity, there exist $s \in \mathbb{N}$ and $\sum_{i=1}^s z_i = 0$ [see the proof of (2.9)], thus $\mu_e(h\mathcal{K}^s\mathbb{1}_A) = 0$. To conclude, we consider the processes $\mathbb{P}_\mu^e$ and $\mathbb{P}_{\mu_e}^e$. Clearly the first is absolutely continuous with respect to the second and the Radon–Nikodym derivative is given by $H((\theta_t)) := h(\theta_0)$. Accordingly, calling $\bar\tau$ the time shift and considering the set $B := \{(\theta_t) \in \Omega : \theta_0 \in A\}$, we have $\mathbb{E}_{\mathbb{P}_{\mu_e}^e}(H\bar\tau^{sn}\mathbb{1}_B) = \mu_e(h\mathcal{K}^{sn}\mathbb{1}_A) = 0$. Thus,

$$0 = \frac{1}{n}\sum_{k=0}^{n-1} \mathbb{E}_{\mathbb{P}_{\mu_e}^e}(H\bar\tau^{sk}\mathbb{1}_B).$$

But the Birkhoff ergodic theorem and the ergodicity of $\mathbb{P}_{\mu_e}^e$ then imply

$$0 = \mathbb{E}_{\mathbb{P}_{\mu_e}^e}(H)\mathbb{E}_{\mathbb{P}_{\mu_e}^e}(\mathbb{1}_B) = \mu_e(A),$$

which is a contradiction. This shows that $\mu$ and $\mu_e$ are equivalent, hence Assumption (A1).

## APPENDIX: CENTRAL LIMIT THEOREM FOR ADDITIVE FUNCTIONALS OF MARKOV CHAINS

In this Appendix, we recall the results and arguments of [19], keeping explicit track of some estimates that are needed for the present paper. The basic idea, going back to [16], is to construct a martingale approximation by solving the equation

$$(1+\varepsilon)h_\varepsilon = \Pi h_\varepsilon + g,$$

where $g = \langle w, G\rangle$ satisfies condition (2.21). Setting $V_t h := \sum_{s=0}^{t-1} \Pi^s h$, a solution to this equation can be written as

$$h_\varepsilon = \varepsilon \sum_{s=1}^{\infty} \frac{V_s g}{(1+\varepsilon)^{s+1}}.$$

If we define

$$H_\varepsilon(\boldsymbol{\omega}, \boldsymbol{\omega}') := h_\varepsilon(\boldsymbol{\omega}') - \Pi h_\varepsilon(\boldsymbol{\omega}),$$

then

$$\sum_{s=0}^{t-1} g(\boldsymbol{\omega}_s) = M_{\varepsilon,t} + \varepsilon \sum_{s=0}^{t-1} h_\varepsilon(\boldsymbol{\omega}_s) + R_{\varepsilon,t},$$

where $M_{\varepsilon,t}$ is an $\mathcal{F}_t$-martingale and $R_{\varepsilon,t}$ a boundary term

$$M_{\varepsilon,t} := \sum_{s=0}^{t-1} H_\varepsilon(\boldsymbol{\omega}_s, \boldsymbol{\omega}_{s+1}); \qquad R_{\varepsilon,t} := \Pi h_\varepsilon(\boldsymbol{\omega}_0) - \Pi h_\varepsilon(\boldsymbol{\omega}_t).$$

LEMMA A.1. *If (2.21) holds, then for each $\varepsilon > 0$, setting $\varepsilon_k := 2^{-k}\varepsilon$, we have $\sum_{k=0}^{\infty} \sqrt{\varepsilon_k}\|h_{\varepsilon_k}\| = o((\ln \varepsilon^{-1})^{-\rho})$, where the norm is the $L^2(\Omega, \mathbb{P}_\mu^e)$-norm.*

PROOF. By definition,

$$\|h_\varepsilon\| \le C + C\varepsilon \sum_{s=2}^{\infty} [\|V_s g\| s^{-3/2} (\ln s)^\rho] \varphi_\varepsilon(s) (\ln s)^{-\rho},$$

where $\varphi_\varepsilon(s) := s^{3/2}(1+\varepsilon)^{-s}$. Then,

$$\sum_{k=0}^{\infty} \sqrt{\varepsilon_k}\|h_{\varepsilon_k}\| \le C\varepsilon^{3/2} + C \sum_{s=2}^{\infty} (\ln s)^{-\rho} [\|V_s g\| s^{-3/2} (\ln s)^\rho] \sum_{k=0}^{\infty} \varepsilon_k^{3/2} \varphi_{\varepsilon_k}(s)$$

and

$$\begin{aligned}\sum_{k=0}^{\infty} \varepsilon_k^{3/2} \varphi_{\varepsilon_k}(s) &\le \sum_{k=0}^{\infty} \left(\frac{s\varepsilon}{2^k}\right)^{3/2} e^{-(s\varepsilon/2^k)(1-(\varepsilon_k/2)e^{\varepsilon_k})} \\ &\le C \int_0^\infty x^{3/2} e^{-x/2} x^{-1}\, dx = C.\end{aligned}$$

Setting $\ell = \varepsilon^{-1/2}$, it follows that

$$\sum_{k=0}^{\infty} \sqrt{\varepsilon}_k \|h_{\varepsilon_k}\| \le C\varepsilon^{3/2} + C\varepsilon^{3/2}\ell^{3/2} \sum_{s=2}^{\ell-1} [\|V_s g\| s^{-3/2} (\ln s)^{\rho}]$$

$$+ C(\ln \ell)^{-\rho} \sum_{s=\ell}^{\infty} [\|V_s g\| s^{-3/2} (\ln s)^{\rho}]$$

$$\le C\varepsilon^{3/4} + C(\ln \varepsilon^{-1})^{-\rho} \sum_{s=\ell}^{\infty} [\|V_s g\| s^{-3/2} (\ln s)^{\rho}]$$

and the result follows trivially from (2.21). □

THEOREM A.2. *If (2.21) holds, then*

$$\sum_{s=0}^{t-1} g(\boldsymbol{\omega}_s) = M_t + R_t,$$

*where* $M_t$ *is a martingale and* $\|R_t\| \le Ct^{1/2}(\ln t)^{-\rho}$.

PROOF. Since, by [19], Lemma 2, $\{H_{\varepsilon_k}\}$ is a Cauchy sequence in $L^2(\mathbf{P})$, we set $H := \lim_{k\to\infty} H_{\varepsilon_k}$. Then, by Lemma A.1, there exist $R_t := \lim_{k\to\infty} R_{\varepsilon_k,t}$ and $\sum_{s=0}^{t-1} g(\boldsymbol{\omega}_s) = M_t + R_t$, where $M_t := \sum_{s=0}^{t-1} H(\boldsymbol{\omega}_s, \boldsymbol{\omega}_{s+1})$ is an $L^2$-martingale. The last estimate follows from Lemma A.1 and equation (8) of [19]. □

**Acknowledgments.** D. D. and C. L. wish to express special thanks to Stefano Olla for many very useful suggestions. We also acknowledge several helpful discussions with Ofer Zeitouni, Carlo Boldrighini and Alessandro Pellegrinotti. We thank Timo Seppalainen for several helpful remarks on a preliminary version of this work. We thank the Institute Henri Poincaré where this work was started (during the trimester *Time at Work*). C. L. thanks the University of Maryland for inviting him during part of this project. Thanks to the Universities of Paris VII and XIII where D. D. and C. L. were visiting during its conclusion.

## REFERENCES

[1] BANDYOPADHYAY, A. and ZEITOUNI, O. (2006). Random walk in dynamic Markovian random environment. *ALEA Lat. Amer. J. Probab. Statist.* **1** 205–224. MR2249655

[2] BOLDRIGHINI, C., MINLOS, R. A. and PELLEGRINOTTI, A. (1994). Central limit theorem for the random walk of one or two particles in a random environment. *Adv. Soviet Math.* **20** 21–75. MR1299595

[3] BOLDRIGHINI, C., MINLOS, R. A. and PELLEGRINOTTI, A. (2000). Random walk in a fluctuating random environment with Markov evolution. In *On Dobrushin's Way. From Probability Theory to Statistical Physics* (R. A. Minlos, S. Shlosman and Yu. M. Suhov, eds.) *Amer. Math. Soc.*, 13–35. Providence, RI. MR1766340

[4] Boldrighini, C., Minlos, R. A. and Pellegrinotti, A. (2004). Random walks in quenched i.i.d. space–time random environment are always a.s. diffusive. *Probab. Theory Related Fields* **129** 133–156. MR2052866

[5] Bolthausen, E. and Sznitman, A.-S. (2002). On the static and dynamic points of view for certain random walks in random environment. *Methods Appl. Anal.* **9** 345–375. MR2023130

[6] Conlon, J. G. and Song, R. (1999). Gaussian limit theorems for diffusion processes and an application. *Stochastic Process. Appl.* **81** 103–128. MR1680489

[7] Comets, F. and Zeitouni, O. (2005). Gaussian fluctuations for random walks in random mixing environments. *Israel J. Math.* **148** 87–113. MR2191225

[8] Derriennic, Y. and Lin, M. (2001). Fractional Poisson equations and ergodic theorems for fractional coboundaries. *Israel J. Math.* **123** 93–130. MR1835290

[9] Derriennic, Y. and Lin, M. (2003). The central limit theorem for Markov chains started at a point. *Probab. Theory Related Fields* **125** 73–76. MR1952457

[10] Fannjiang, A. and Komorowski, T. (2001). Invariance principle for a diffusion in a Markov field. *Bull. Polish Acad. Sci. Math.* **49** 45–65. MR1824156

[11] Gordin, M. I. (1969). The central limit theorem for stationary processes. *Dokl. Akad. Nauk SSSR* **188** 739–741. MR0251785

[12] Grimmett, G. and Stirzaker, D. (2001). *Probability and Random Processes*, 3rd ed. Oxford Univ. Press. MR2059709

[13] Hall, P. and Heyde, C. C. (1980). *Martingale Limit Theory and Its Applications*. Academic Press, New York. MR0624435

[14] Komorowski, T. and Olla, S. (2001). On homogenization of time-dependent random flows. *Probab. Theory Related Fields* **121** 98–116. MR1857110

[15] Keller, G. and Liverani, C. (2006). Uniqueness of the SRB measure for piecewise expanding weakly coupled map lattices in any dimension. *Comm. Math. Phys.* **262** 33–50. MR2200881

[16] Kipnis, C. and Varadhan, S. R. S. (1986). Central limit theorem for additive functionals of reversible Markov processes and applications to simple exclusions. *Comm. Math. Phys.* **104** 1–19. MR0834478

[17] Landim, C., Olla, S. and Yau, H. T. (1998). Convection–diffusion equation with space–time ergodic random flow. *Probab. Theory Related Fields* **112** 203–220. MR1653837

[18] Liverani, C. (1996). Central limit theorem for deterministic systems. In *International Conference on Dynamical Systems* (*Montevideo, 1995*). *A Tribute to Ricardo Mane.* (F. Ledrappier, J. Levovicz and S. Newhouse, eds.) 56–75. Longman, Harlow. MR1460797

[19] Maxwell, M. and Woodroofe, M. (2000). Central limit theorems for additive functionals of Markov chains. *Ann. Probab.* **28** 713–724. MR1782272

[20] Rassoul-Agha, F. and Seppalainen, T. (2005). An almost sure invariance principle for random walks in a space–time i.i.d. random environment. *Probab. Theory Related Fields* **133** 299–314. MR2198014

[21] Rassoul-Agha, F. and Seppalainen, T. (2005). An almost sure invariance principle for additive functionals of Markov chains. Preprint. Available at http://arxiv.org/abs/math/0411603v2.

[22] Stannat, W. (2004). A remark on the CLT for random walks in a random environment. *Probab. Theory Related Fields* **130** 377–387. MR2095935

[23] Sinai, Y. G. (1993). A random walk with a random potential. *Teor. Veroyatnost. i Primenen.* **38** 457–460. [Translation in *Theory Probab. Appl.* **38** (1993) 382–385.] MR1317991

[24] SINAI, Y. G. (1995). A remark concerning random walks with random potentials. *Fund. Math.* **147** 173–180. MR1341729

[25] SZNITMAN, A.-S. (2004). Topics in random walk in random environment. In *School and Conference on Probability Theory* 203–266. Electronic. MR2198849

[26] ZEITOUNI, O. (2004). Random walks in random environment. *Lectures on Probability Theory and Statistics. Lecture Notes in Math.* **1837** 189–312. Springer, Berlin. MR2071631

D. DOLGOPYAT
DEPARTMENT OF MATHEMATICS
UNIVERSITY OF MARYLAND
4417 MATHEMATICS BLDG
COLLEGE PARK, MARYLAND 20742
USA
E-MAIL: dmitry@math.umd.edu

G. KELLER
MATHEMATISCHES INSTITUT
UNIVERSITÄT ERLANGEN-NÜRNBERG
BISMARCKSTR. 1 1/2
91054 ERLANGEN
GERMANY
E-MAIL: keller@mi.uni-erlangen.de

C. LIVERANI
DIPARTIMENTO DI MATEMATICA
UNIVERSITÀ DI ROMA (TOR VERGATA)
VIA DELLA RICERCA SCIENTIFICA
00133 ROMA
ITALY
E-MAIL: liverani@mat.uniroma2.it